\documentclass[12pt]{amsart}
\usepackage{amsmath,amsfonts,amssymb,amscd,amsthm,amsbsy,epsf}
\newtheorem{thm}{Theorem}[section]
\newtheorem{lem}[thm]{Lemma}
\newtheorem{cor}[thm]{Corollary}
\newtheorem{prp}[thm]{Proposition}
\newtheorem*{Hindman}{Hindman's Theorem}    
\newtheorem*{HindmanSet}{Hindman's Theorem (Set Version)}   
\newtheorem*{HindMill}{Hindman-Milliken Theorem (Set Version)}  
\theoremstyle{definition}
\newtheorem{exa}[thm]{Example}
\newtheorem{fct}[thm]{Fact}
\newtheorem{rmk}[thm]{Remark}

\newtheorem {Obs}[thm]{Observations}

\renewcommand \qedsymbol {$\boldsymbol{\dashv}$} 

\def\real{{\mathbb R}}
    \def\Norm{{|\!|\!|}}
    \def\NORM{{|\!|\!|\!|}}
\def\BigNorm{{\Big|\!\Big|\!\Big|}}
    
\def\eqref#1{{\rm (\ref{#1})}}

\newcommand {\art}[6]{{\sc #1:} {\it #2}, {\bfseries\itshape #3}, {vol.\;#4} {(#5),} {#6.}}
\newcommand {\bok}[2]{{\sc #1:} {\bfseries\itshape #2}, {\sf{to appear.}}}
\newcommand {\book}[5]{{\sc #1:} {\bfseries\itshape #2}, {#3,} {#4} {(#5).}}
\newcommand {\notes}[5]{{\sc #1:} {\it #2}, {\bfseries\itshape #3}, {#4,} {#5.}}
\newcommand {\samp}[8]{{\sc #1:} {\it #2}, {\bfseries\itshape #3}, {(#4, {\sl Ed.}),}
                                  {#5,} {#6,} {(#7),} {#8.}}
\newcommand {\preprint}[2]{{\sc #1:} {\it #2}, {\sf{preprint.}}}

\newcommand \spart {\la\omega\ra^\omega}
\newcommand \reals {[\omega ]^\omega}
\newcommand \fin {[\omega]^{<\omega}}

\newcommand \te {{\text{\sl th}}}
\newcommand \jump {{\hfill\\[0.6ex]}}

\newcommand \ind {\hspace*{0.9ex}}
\newcommand \mind {\hspace*{-0.9ex}}
\newcommand \indd {\hspace*{1.4ex}}
\newcommand \mindd {\hspace*{-1.4ex}}
\newcommand \inddd {\hspace*{2.9ex}}
\newcommand \minddd {\hspace*{-2.9ex}}

\newcommand \supp {\operatorname{supp}}

\newcommand \ceq {\sqsubseteq}
\newcommand \afeq {{}^*\hspace*{-0.8ex}\sqsupseteq}

\newcommand \eps {\varepsilon}

\newcommand \la {\langle}
\newcommand \ra {\rangle}
\newcommand \subs {\subseteq}

\newcommand \cF {{\mathcal {F}}}
\newcommand \cG {{\mathcal {G}}}

\newcommand \bfa {\boldsymbol{a}}

\newtheorem {subnummer}{ }[thm]
\newtheorem {subnr}[subnummer]{{\rm\small{(}}\!\rm\small}
\newcommand {\beginsubnr}[1]
            {\begin{subnr}\label{#1}\hspace*{-1ex}{\rm\small{)}}\hspace*{1ex}}

\begin{document}
\begin{center}
      {\large{\bf{ON ASYMPTOTIC MODELS IN BANACH SPACES}}}
\end{center}
\medskip

\begin{center}
{\sc{Lorenz Halbeisen}}\\[1.1ex] {\small{\sl Dept. of Pure
Mathematics\\ Queen's University Belfast\\ Belfast BT7 1NN (Northern
Ireland)}}\\[0.8ex] {\small{\sl Email: halbeis@qub.ac.uk}}
\end{center}
\medskip

\begin{center}
{\sc{Edward Odell}\footnote{Research supported by NSF.}}\\[1.1ex]
{\small{\sl Department of Mathematics\\ University of Texas\\ Austin,
TX 78712-1082 (U.S.A.)}}\\[0.8ex] {\small{\sl Email:
odell@math.utexas.edu}}
\end{center}
\medskip

\setlength{\rightskip}{1cm}  \setlength{\leftskip}{1cm} {\small {\bf
Keywords:} Ramsey's Theorem, Hindman-Milliken Theorem, block
sequences, spreading models, asymptotic structure, asymptotic
models.} \hfill\smallskip

{\small {\bf 2000 Mathematics Subject Classification:} {\bf 46B45
05D10} 46B35 05D05 46B20} \hfill\bigskip


\section*{Abstract}

A well known application of Ramsey's Theorem to Banach Space Theory
is the notion of a spreading model $(\tilde{e}_i)$ of a normalized
basic sequence $(x_i)$ in a Banach space $X$. We show how to
generalize the construction to define a new creature $(e_i)$, which
we call an asymptotic model of $X$. Every spreading model of $X$ is
an asymptotic model of $X$ and in most settings, such as if $X$ is
reflexive, every normalized block basis of an asymptotic model is
itself an asymptotic model. We also show how to use the
Hindman-Milliken Theorem---a strengthened form of Ramsey's
Theorem---to generate asymptotic models with a stronger form of
convergence.

\setcounter{section}{-1}

\section{Introduction}\label{sec:introduction}

Ramsey Theory, and especially Ramsey's Theorem, is a very powerful
tool in infinitary combinatorics and has many interesting (and
sometimes unexpected) applications in various fields of Mathematics.
Generally speaking, theorems in Ramsey Theory are of the type that a
function into a finite set can be restricted to some sort of infinite
substructure, on which it is constant. In applications to analysis we
successively apply Ramsey's Theorem to certain $\eps$-nets to obtain
infinite substructures on which certain Lipschitz functions are
nearly constant in an asymptotic sense (cf. e.g., \cite{Odell} or
\cite[Part\,III]{OdellKechris}).

A well known application of Ramsey's Theorem
(\cite[Theorem\;A]{Ramsey}) to Banach Space Theory is due to
Antoine~Brunel and Louis~Sucheston (cf. \cite{BrunelSucheston}).
Roughly speaking, it says that every normalized basic sequence in a
Banach space has a subsequence which is ``asymptotically''
subsymmetric, ultimately yielding a spreading model.

There are two main directions to generalize Ramsey's Theorem. One is
in terms of partitions and another one leads to the so-called Ramsey
property. (Some results concerning the symmetries between the
combination of these two directions can be found in \cite{sym}.) Both
directions are already used in Banach Space Theory. For example the
fact that Borel sets have the Ramsey property is used in Farahat's
proof of Rosenthal's Theorem, which says that a normalized sequence
has a subsequence which is either equivalent to the unit vector basis
of $\ell_1$ or is weakly Cauchy. Further, a combination of both
directions is used by Gowers in the proof of his famous
Dichotomy~Theorem.

In the sequel, we prove a generalized version of the Brunel-Sucheston
Theorem by using Ramsey's Theorem.
We apply this to basic arrays, namely
certain sequences of basic sequences in
$X$. Also we show how a generalization of Ramsey's Theorem, the
Hindman-Milliken Theorem, can be used to construct asymptotic models with
a stronger form of convergence.

The object we obtain, a basis $(e_i)_{i\in\omega}$ for some infinite
dimensional Banach space $E$, we call an asymptotic model of $X$.
Asymptotic models include not only all spreading models of $X$, and
even in many cases all normalized block bases of such, but more general
sequences as well.
If the sequences in the generating basic array are all block bases
of a fixed basis or are all weakly null then
the notion lies somewhere between that of spreading
models and asymptotic structure (see\;\cite{Maurey.etal}), although it is
closer in flavor to the theory of spreading models.
The construction we use to get an asymptotic model has been used in the
past by several authors to study spreading models and the behavior of
sequences over $X$ (e.g., \cite{Ro83}, \cite{Ma83} and
\cite{Androulakis.etal}.
In particular in \cite{Ro83} the concept of an $\infty$-type over a
Banach space is introduced and this actually contains within it the
notion  of an asymptotic model.
But our more restricted viewpoint in this paper is the first study
of what we have chosen to call ``asymptotic models'' themselves.

In Section\;\ref{sec:HindmanMilliken} we recall the Hindman-Milliken
Theorem. In Section\;\ref{sec:models} we define and construct
asymptotic models. In addition we make a number of observations about
asymptotic models and their relation with spreading models and
asymptotic structure. Section\;\ref{sec:renorm} generalizes some
results of
\cite{OS98b} 
to the setting of asymptotic models.
Section\;\ref{sec:odd} concerns some stronger
versions one might hope to have, but as we show one cannot achieve in
general. In this section we also raise some open problems.

For the reader's convenience, we recall some set theoretic
terminology we will use frequently. A natural number $n$ is
considered as the set of all natural numbers less than $n$, in
particular, $0=\emptyset$. Let $\omega =\{0,1,2,\ldots\}$ denote the
set of all natural numbers. By the way, we always start counting by
0. Some more set theoretic terminology will be introduced in the
following section.

The notation concerning sequence spaces is standard and can be found
in textbooks like \cite{Diestel}, \cite{Guerre} and
\cite{LindenstraussTzafriri}. However, for the sake of the
non-expert, we recall some definitions.

A sequence $(x_i)_{i\in\omega}$ in a normed space is {\bf normalized}
if for all $i\in\omega$, $\|x_i\|=1$, and it is {\bf seminormalized}
if there exists an $M$ with $0<M<\infty$ such that all $i\in\omega$,
$\frac{1}{M}\le \|x_i\|\le M$. If $(x_i)_{i\in\omega}$ is a sequence
of non-zero vectors in a Banach space $X$, then $(x_i)_{i\in\omega}$
is {\bf basic} iff there exists $C<\infty$ so that for all $n<m$ and
$(a_i)_{i\in m} \subseteq \real$, $\|\sum_{i\in n} a_i x_i\| \le C\|
\sum_{i\in m} a_ix_i\|$. The smallest such $C$ is called the {\bf
basis constant} of $(x_i)_{i\in\omega}$ and $(x_i)_{i\in\omega}$ is
then called {\bf $\boldsymbol{C}$-basic}. The basic sequence
$(x_i)_{i\in\omega}$ is {\bf monotone basic\/} if it is $1$-basic,
and it is {\bf bimonotone\/} if it is monotone and the tail
projections are monotone as well (i.e., $I-P_n$ has norm one if $P_n$
it the $n^\te$ initial projection). If $(x_i)_{i\in\omega}$ is basic,
then every $x$ in the closed linear span of $(x_i)_{i\in\omega}$ can
be uniquely expressed as $\sum_{i\in\omega} a_i x_i$ for some
$(a_i)_{i\in\omega}\subseteq \real$. Basic sequences
$(x_i)_{i\in\omega}$ and $(y_i)_{i\in\omega}$ are {\bf
$\boldsymbol{C}$-equivalent\/} if there exist constants $A$ and $B$
with $AB\le C$ so that for all $n\in\omega$ and scalars $(a_i)_{i\in
n}$ $$A^{-1}\Big{\|}\sum_{i\in n} a_i x_i \Big{\|}\le
\Big{\|}\sum_{i\in n} a_i y_i \Big{\|}\le B\Big{\|}\sum_{i\in n} a_i
x_i \Big{\|}\,.$$ For a basic sequence $(x_i)_{i\in\omega}$ and
scalars $(b_l)_{l\in\omega}$, a sequence of non-zero vectors
$(y_j)_{j\in\omega}$ of the form $$y_j={\sum_{l=p_k}^{p_{k+1}-1}b_l
x_l}\,,$$ where $p_0<p_1<\ldots<p_k<\ldots$ is an increasing sequence
of natural numbers, is called a {\bf block basic sequence\/} or just
a {\bf block basis}.

A basic sequence $(x_i)_{i\in\omega}$ is called {\bf boundedly
complete} if, for every sequence of scalars $(a_i)_{i\in\omega}$ such
that $\sup_n \big{\|} \sum_{i\in n}a_i x_i\big{\|}<\infty$, the
series $\sum_{i\in\omega} a_i x_i$ converges. A basic sequence
$(x_i)_{i\in\omega}$ is {\bf unconditional\/} if for any sequence
$(a_i)_{i\in\omega}$ of scalars and for any permutation $\pi$ of
$\omega$, i.e., for any bijection $\pi:\omega\to\omega$,
$\sum_{i\in\omega}a_i x_i$ converges if and only if
$\sum_{i\in\omega}a_{\pi(i)} x_{\pi(i)}$ converges. A non-zero
sequence of vectors $(x_i)_{i\in\omega}$ is {\bf unconditional
basic\/} iff there exists $C<\infty$ so that for all $n\in\omega$,
$\eps_i =\pm 1$ and $(a_i)_{i\in n} \subseteq \real$, $\|\sum_{i\in
n} \eps_i a_i x_i\| \le  C\|\sum_{i\in n} a_i x_i\|$. The smallest
such $C$ is the {\bf unconditional basis constant\/} of $(x_i)$.

A normalized basic sequence $(x_i)_{i\in\omega}$ is {\bf
$\boldsymbol{C}$-subsymmetric} if $(x_i)_{i\in\omega}$ is
$C$-equivalent to each of its subsequences (notice that we do not
require it to be unconditional which differs from the terminology of
\cite{LindenstraussTzafriri}).

For a set of vectors $A$, $\langle A\rangle$ denotes the linear span
of $A$ and $\big[ A \big]$ denotes the {\it closure of the linear
span\/} of $A$. Note that if the normalized basic sequences
$(x_i)_{i\in\omega}$ and $(y_i)_{i\in\omega}$ are $C$-equivalent,
then the spaces $[(x_i)_{i\in\omega}]$ and $[(y_i)_{i\in\omega}]$ are
$C$-isomorphic.

The {\it dual space\/} of a Banach space $X$ is denoted by $X^*$.

Suppose that $(x_i)_{i\in\omega}$ is a basic sequence. For each $x^*$
in $[(x_i)_{i\in\omega}]^*$ and each $n\in\omega$, let
$\|x^*\|_{(n)}$ be the norm of the restriction of $x^*$ to $\big{[}\{
x_i:i>n\}\big{]}$. Then $(x_i)_{i\in\omega}$ is {\bf shrinking\/} if
for each $x^*\in [(x_i)_{i\in\omega}]^*$,
$\lim_{n\to\infty}\|x^*\|_{(n)}=0$.

If $Y$ is a normed linear space, $B_Y$ denotes the {\it closed unit
ball\/} of $Y$ and $S_Y$ is the {\it unit sphere}. In the sequel, $X$
will always denote a separable infinite dimensional real Banach
space.

\section{Special Partitions}
\setcounter{equation}{0}
\setcounter{thm}{0}

Let $\omega +1:=\omega\cup\{\omega\}$, so if $\eta \in\omega +1$,
then $\eta$ is either a natural number or $\eta = \omega$. If $x$ is
a set, we write $|x|$ for the cardinality of $x$. We will use
$\omega$ also as a cardinal number, namely $\omega=|\omega |$. If $x$
is a set and $\eta \in\omega+1$, then $$[x]^\eta
:=\{y\subs x:|y|=\eta\}$$ and $$[x]^{<\eta} :=\{y\subs
x:|y|<\eta\}\,.$$

If $a,b\subs\omega$, we write $a<b$ in place of ``for all $n\in a$
and $m\in b$, $n<m$''. Note that $a<b$ implies $a\in\fin$.

A {\bf partition\/} $P$ of set $S$ is a set of non-empty, pairwise
disjoint subsets of $S$ such that $\bigcup P=S$. For a partition $P$,
the sets $b\in P$ are called the {\bf blocks\/} of $P$.

In the following we consider ``special'' partitions of subsets of
$\omega$.

If $P$ is a partition of some subset of $\omega$, then $P$ is called
a {\bf special partition\/}, if for all blocks $a,b\in P$ we have
either $a<b$, or $a=b$, or $a>b$.

Notice that if $P$ is a special partition with infinitely many
blocks, then all of its blocks are finite.

For $\eta\in\omega +1$, let $\la\omega\ra^\eta$ denote the set of all
special partitions of subsets of $\omega$ such that $|P|=\eta$. In
particular, $\spart$ is the set of all special partitions with
infinitely many blocks.

Let $P_1,P_2$ be two special partitions. We say that $P_1$ is {\bf
coarser} than $P_2$, or that $P_2$ is {\bf finer\/} than $P_1$, and
write $P_1\ceq P_2$, if each block of $P_1$ is the union of blocks of
$P_2$.

For a special partition $P$ and $\eta\in\omega+1$ let $$\la
P\ra^\eta:=\{Q: Q\ceq P\wedge |Q|=\eta\}\,.$$

If $P$ is a special partition and $b\in P$, then $\min(b):=\bigcap b$
denotes the minimum of the set $b$. If we order the blocks of $P$ by
their minimum, then $P(n)$ denotes the $n$th block with respect to
this ordering.

If $P_1,P_2$ are two special partitions, then we write $P_1\ceq^*
P_2$ if there is an $n\in\omega$ such that
$$\big(P_1\setminus\{P_1(i):i\in n\}\big)\ceq P_2\,.$$
In other
words, $P_1\ceq^* P_2$ if all but finitely many blocks of $P_1$ are
unions of blocks of $P_2$.

\begin{fct}\label{fct:tower} If $P_0\,\afeq P_1\,\afeq P_2\,\afeq
\ldots\afeq P_i\,\afeq\ldots$ where $P_i\in\spart$ (for each
$i\in\omega$), then there is a special partition $P\in\spart$ such
that for each $i\in\omega$, $P\ceq^* P_i$.
\end{fct}

(The proof is similar to the proof of Fact\;2.3 of \cite{sym}.)

\section{The Hindman-Milliken Theorem}\label{sec:HindmanMilliken}
\setcounter{equation}{0}
\setcounter{thm}{0}

First, we recall the well-known Hindman Theorem, and then we give
Milliken's generalization of Hindman's Theorem.

If $A\in\fin$, then we write $\sum A$ for $\sum_{a\in A} a$, where we
define $\sum\emptyset:=0$.

In \cite{Hindman}, Neil~Hindman proved the following.

\begin{Hindman}
If $m$ is a positive natural number and
$f:\omega\to m$ is a function, then there exist $r\in m$ and
$x\in\reals$ such that whenever $A\in [x]^{<\omega}$ is non-empty, we
have $f (\sum A) =r$.
\end{Hindman}

Ronald~Graham and Bruce~Rothschild noted that Hindman's Theorem can
be formulated in terms of finite sets and their unions instead of
natural numbers and their sums. This yields the following.

\begin{HindmanSet}
If $m$ is a positive natural
number, $I\in\reals$ and $f:[I]^{<\omega}\to m$ is a function, then
there exist $r\in m$ and an infinite set $H\subs [I]^{<\omega}$ such
that $a\cap b=\emptyset$ for all distinct sets $a,b\in H$, and
whenever $A\in [H]^{<\omega}$ is non-empty, we have $f (\bigcup A) =r$.
\end{HindmanSet}

Using Hindman's Theorem as a strong pigeonhole principle,
Keith~Milliken proved a strengthened version of Ramsey's Theorem,
which we will call the Hindman-Milliken Theorem (cf.
\cite[Theorem\;2.2]{Milliken}). The Hindman-Milliken Theorem in terms
of unions can be stated as follows:

\begin{HindMill}
Let $m,n$ be positive
natural numbers, $Q\in\spart$ and $f:\la Q\ra^n\to m$ a function,
then there is an $P\in\la Q\ra^\omega$ such that $f$ is constant on
$\la P\ra^n$.
\end{HindMill}

As consequences of the Hindman-Milliken Theorem one gets Ramsey's
Theorem (Theorem\;A of \cite{Ramsey}) as well as Hindman's Theorem
(cf. \cite{Milliken}).

\section{Asymptotic Models}\label{sec:models}
\setcounter{equation}{0}
\setcounter{thm}{0}

First we recall the notion of a spreading model. If
$(x_i)_{i\in\omega}$ is a normalized basic sequence in a Banach space
$X$ and $\eps_n\downarrow 0$ (a sequence of positive real numbers
which tends to $0$), then one can find a subsequence
$(y_i)_{i\in\omega}$ of $(x_i)_{i\in\omega}$ such that the following holds:
For any positive $n\in\omega$, any sequence $(a_k)_{k\in n}\in
[-1,1]^{n}$ and any natural numbers $n\le i_0<\ldots
<i_{n-1}$ and $n\le j_0<\ldots <j_{n-1}$ we have
$$\bigg|\Big{\|}
\sum_{k\in n} a_k y_{i_k}\Big{\|} -\Big\| \sum_{k\in n} a_k
y_{j_k}\Big{\|}\bigg| < \eps_n\ .$$

This is proved by using Ramsey's Theorem iteratively for a finite
$\delta_n$-net in the unit ball of $\ell_\infty^n$ ($\delta_n$
depends upon $\eps_n$) to stabilize, up to $\delta_n$, the functions
$f(i_0,\ldots,i_{n-1})\equiv \|\sum_{i\in n} a_i x_i\|$ over a
subsequence $(y_i)_{i\in\omega}$ of $(x_i)_{i\in\omega}$ for each
$(a_i)_{i\in n}$ in the $\delta_n$-net. Thus, one obtains a limit,
$\|\sum_{i\in n} a_i \tilde{e}_i\|$, for each finite sequence
$(a_i)_{i\in n}$ of scalars. The sequence
$(\tilde{e}_i)_{i\in\omega}$ is called a {\bf spreading model of\/}
$(y_i)_{i\in\omega}$; $(\tilde{e}_i)_{i\in\omega}$ is a normalized
1-subsymmetric basis for $\tilde{E}$, the closed linear span of the
$\tilde{e}_i$'s, and $\tilde{E}$ is called a {\bf spreading model\/}
of $X$ generated by $(\tilde{e}_i)_{i\in\omega}$. Hence, for any
natural numbers $j_0<\ldots
<j_{n-1}$ we have $\|\sum_{i\in n} a_i \tilde{e}_i\|
=\|\sum_{i\in n} a_i \tilde{e}_{j_i}\|$. If
$(y_i)_{i\in\omega}$ is weakly null, $(\tilde{e}_i)_{i\in\omega}$ is
suppression-$1$ unconditional: $\|\sum_{i\in F} a_i
\tilde{e}_i\| \le \|\sum_{i\in\omega} a_i
\tilde{e}_i\|$ for all $F\subs\omega$ and each sequence
$(a_i)_{i\in\omega}$ of scalars. These facts can be found in
\cite{BeauzamyLapreste} or \cite{Odell}.

Before presenting our extension we set some notation.

We shall call $(x_i^n)_{n,i\in\omega}$ a {\bf $\boldsymbol{K}$-basic
array\/} in $X$, if for all $n\in\omega$, $(x_i^n)_{i\in\omega}$ is a
$K$-basic normalized sequence in $X$ and moreover if for all
$m\in\omega$ and all integers $m\le i_0<\ldots <i_{m-1}$, every
sequence  $(x_{i_j}^j)_{j\in m}$ is $K$-basic.
Furthermore,
$(x_i^n)_{n,i\in\omega}$ is a {\bf basic array\/} in $X$ if it is a
$K$-basic array for some $K<\infty$.

If $X$ has a basis $(x_i)_{i\in\omega}$ then $(x_i^n)_{n,i\in\omega}$
is a {\bf block basic array\/} in $X$ (with respect to
$(x_i)_{i\in\omega}$) if in addition each row $(x_i^n)_{i\in\omega}$
is a block basis of $(x_i)_{i\in\omega}$ and all sequences
$(x_{i_j}^j)_{j\in m}$ as described above are also block bases of
$(x_i)_{i\in\omega}$.

In what we present, the only important part of the array is the upper
triangular part: $\{x_i^n: n\in\omega$ and $i\ge n\}$. The lower
triangular part can be ignored or omitted and we shall often do so.

\begin{prp}\label{prp:K-basic}
Let $(x_i^n)_{n,i\in\omega}$ be a $K$-basic array in some Banach
space $X$. Then given $\eps_n\downarrow 0$, there exists a
subsequence $(k_n)_{n\in\omega}$ of $\omega$ so that for all
$n\in\omega$, $(b_i)_{i\in n}\in [-1,1]^n$, $n\le i_0<\cdots <
i_{n-1}$ and $n\le \ell_0 < \cdots < \ell_{n-1}$,
$$\left| \bigg\|
\sum_{j\in n} b_j x_{k_{i_j}}^j\bigg\| - \bigg\| \sum_{j\in n} b_j
x_{k_{\ell_j}}^j \bigg\| \right| < \eps_n\ .$$
\end{prp}

\begin{proof}
As in the case of spreading models, this follows easily from Ramsey's Theorem
and the standard diagonalization argument.
One $\frac{\eps_n}2$-stabilizes $f(i_0,\ldots,i_{n-1}) := \|\sum_{j\in n}
b_j x_{i_j}^j\|$
over all subsequences of length $n$ on some subsequence of $\omega$ for
each of finitely many $(b_j)_{j\in n} \in [-1,1]^n$ out of some
$\delta_n$-net in $B_{\ell_\infty^n}$.
\end{proof}

If the conclusion of the proposition holds for
$(y_i^n)_{n,i\in\omega}$, where $y_i^n = x_{k_i}^n$, then the
iterated limit, $\lim_{i_0\to\infty}\cdots \lim_{i_{n-1}\to\infty}
\|\sum_{j\in n} b_j  y_{i_j}^j\|$, defines a norm on $c_{00}$, the
linear space of finitely supported real sequences on $\omega$. We let
$E$ be the completion of $c_{00}$ under this norm. The unit vector
basis $(e_i)_{i\in\omega}$ thus becomes a $K$-basis for $E$. We call
$(e_i)_{i\in\omega}$ or $E$ an {\bf asymptotic model\/} of $X$ {\bf
generated by\/} $(y_i^n)_{n,i\in\omega}$.

If $(x_i^n)_{n,i\in\omega}$ is a basic array and $i_0<i_1<\cdots$
then $(y_j^n)_{n,j\in\omega}$, where $y_j^n = x_{i_j}^n$, is called a
{\bf subarray\/} of $(x_i^n)_{n,i\in\omega}$.
Proposition\;\ref{prp:K-basic} says that every basic array admits a
subarray which generates an asymptotic model. Also clearly if
$(y_i^n)_{n,i\in\omega}$ generates $(e_i)_{i\in\omega}$, then every
subarray of $(y_i^n)_{n,i\in\omega}$ generates $(e_i)$ as well.

We shall have occasion to use the following simple lemma.

\begin{lem}\label{lem:normalized}
For each $n\in\omega$ let $(x_i^n)_{i\in\omega}$ be a normalized
sequence in a Banach space $X$. If either
\begin{itemize}
\item[a)] each $(x_i^n)_{i\in\omega}$ is weakly null or
\item[b)] each $(x_i^n)_{i\in\omega}$ is a block basis of some basic
sequence $(x_i)_{i\in\omega}$ in $X$,
\end{itemize}
then the array $(x_i^n)_{n,i\in\omega}$ admits a basic subarray
$(y_i^n)_{n,i\in\omega}$. {If\/} {\bf a)}, then given $\eps>0$,
$(y_i^n)_{n,i\in\omega}$ can be chosen to be a $1+\eps$-basic array.
{If\/} {\bf b)}, $(y_i^n)_{n,i\in\omega}$ can be chosen to be a block
basic array of $(x_i)_{i\in\omega}$.
\end{lem}

\begin{proof}
To prove b) we need just choose the subarray $(y_i^n)_{n,i\in\omega}$
so that for all $n\in\omega$, $j\in n$, $i\in n+1$, $\max \big(\supp
(y_{n-1}^j)\big) < \min \big(\supp (y_n^i)\big)$ where if $y= \sum
a_i x_i$ then $\supp (y) = \{i: a_i\ne 0\}$. a)~is proved by a slight
generalization of the proof of the well known fact that a normalized weakly
null sequence admits a $1+\eps$-basic subsequence.
One takes $\eps_n\downarrow 0$ rapidly and then chooses the column
$(y_n^i)_{i\in\omega}$ so that $|f(y_n^i)| <\eps_n$ for $i\in n+1$ and
each $f$ in a finite $1+\eps_n$-norming set of functionals of $B_
{\langle y_i^j : i,j\in n\rangle}$
\end{proof}

We will call a basic array $(x_i^n)$ whose rows,
$(x_i^n)_{i\in\omega}$, are all weakly null a {\bf weakly null basic
array}.

If $(e_i)_{i\in\omega}$ is a spreading model of $X$ generated by the
basic sequence $(x_i)_{i\in\omega}$, then clearly
$(e_i)_{i\in\omega}$ is an asymptotic model of $X$ as well (generated
by $(x_i^n)_{n,i\in\omega}$ where $x_i^n = x_i$ for all
$n,i\in\omega$). A block basis of a spreading model need not be a
spreading model, however, this is not usually the case for asymptotic
models. But first  we introduce some new notation and a new stronger
way of obtaining asymptotic models.

A basic array is a {\bf strong $\boldsymbol{K}$-basic array\/} if in
addition to the defining conditions of a $K$-basic array, for all
integers $m\le i_0 < i_1 <\cdots < i_{m-1}$, every sequence of
non-zero vectors $(y_j)_{j\in m}$ is $K$-basic whenever $y_j \in
\langle x_s^j :i_j \le s<i_{j+1}\rangle$. Note that the proof of
Lemma\;\ref{lem:normalized} actually yields that one can choose the
subarray $(y_i^n)_{n,i\in\omega}$ to be strong basic.

Let $(x_i^n)_{n,i\in\omega}$ be a strong basic array. Given
$m\in\omega$, a finite set of positive integers
$F=(i_0,i_1,\ldots,i_{n-1})$ with $i_0<\ldots
<i_{n-1}$, and a (possibly infinite) sequence
${\bfa}=(a_0,a_1,\ldots)$ of scalars of length at least $n$ with
$a_i\neq 0$ for some $i\in n$, we define
$$x^m(F,{\bfa}):=\frac{\sum_{j\in n} a_j x_{i_j}^m}{\big\| \sum_{j\in
n} a_j x_{i_j}^m\big\|}\ .$$

\begin{thm}\label{thm:main}
Let $X$ be a Banach space and let $(x_i^n)_{n,i\in\omega}$ be a
strong $K$-basic array in $X$ for some $K<\infty$.
For $i\in\omega$ and each non-empty finite
set of integers $F=\{i_0,\ldots,i_{n-1}\}$ with $i_0<\ldots
<i_{n-1}$, let ${\bfa}_F^i$ be a (possibly infinite) sequence of scalars
of length at least $n$ and not identically zero in the first $n$
coordinates and let $\eps_n \downarrow 0$. Then there exists a
special partition $P=\{P(k):k\in\omega\}\in\spart$ such that the
following holds: For all positive $n\in\omega$ and $(b_i)_{i\in n}\in
[-1,1]^{n}$ and $s,t\in\la P\ra^n$ with $\min (s(0)), \min (t(0))\ge
n$ we have $$\bigg|\Big\| \sum_{i\in n}
b_i\,x^i\big(s(i),{\bfa}_{s(i)}^i\big)\Big\| -\Big\| \sum_{i\in n}
b_i\,x^i\big(t(i),{\bfa}_{t(i)}^i\big)\Big\|\bigg| < \eps_n\ .$$
\end{thm}

\begin{proof}
The theorem follows from the Hindman-Milliken Theorem the same way
that one obtains a subsequence of a given basic sequence
$(x_i)_{i\in\omega}$ yielding a spreading model via Ramsey's Theorem:
Given finitely many sequences $(b_i)_{i\in n}\in [-1,1]^n$, a
$\delta_n$-net in $B_{\ell^n_\infty}$ (the unit ball of
$\ell_{\infty}^n$) for an appropriate $\delta_n$, and a special
partition $P\in\spart$, then one can find $Q\in\la P\ra^\omega$ so
that for all $t,r\in\la Q\ra^n$ we have
\begin{equation}
\bigg| \Big\| \sum_{i\in n} b_i\,x^i\big(t(i),
{\bfa}_{t(i)}^i\big)\Big\| -\Big\| \sum_{i\in n}
b_i\,x^i\big(r(i), {\bfa}_{r(i)}^i\big)\Big\| \bigg| <
\delta_n\ .\tag{$*$}
\end{equation}
One then uses standard approximation and diagonalization (see
Fact\;\ref{fct:tower}) arguments to conclude the proof.\hfill\\
Indeed, given $(b_i)_{i\in n}$ and a special partition $P\in\spart$,
we partition the interval $[-n,n]$ into say $m$ disjoint subintervals
$(I_i)_{i\in m}$, each of length less than $\delta_n$. Given $t\in\la
P\ra^n$, we let $$f(t):=j\ \,\text{if and only if}\ \,\Big\|
\sum_{i\in n} b_i\,x^i\big(t(i),{\bfa}_{t(i)}^i\big)\Big\| \in I_j\
.$$ An application of the Hindman-Milliken Theorem yields $Q \in\la
P\ra^\omega$ so that $(*)$ holds for all $t,r\in\la Q\ra^n$. We
repeat this for each $(b_i)_{i\in n}$. For an arbitrary $(c_i)_{i\in
n}\in [-1,1]^n$ one chooses $(b_i)_{i\in n}$ from this $\delta_n$-net
with $|c_i-b_i|<\delta_n$ (for all $i\in n$). Hence, for $t,r\in\la Q
\ra^n$,
 \begin{multline*}
\bigg| \Big\| \sum_{i\in n} c_i\,x^i\big(t(i),{\bfa}_{t(i)}^i\big)\Big\| -
 \Big\| \sum_{i\in n} c_i\,x^i\big(r(i),{\bfa}_{r(i)}^i\big)\Big\|
 \bigg| =\\
 \bigg| \Big\| \sum_{i\in n} c_i\,x^i\big(t(i),{\bfa}_{t(i)}^i\big) -
          \sum_{i\in n} b_i\,x^i\big(t(i),{\bfa}_{t(i)}^i\big) +
          \sum_{i\in n} b_i\,x^i\big(t(i),{\bfa}_{t(i)}^i\big)\Big\| \\
-\Big\| \sum_{i\in n} c_i\,x^i\big(r(i),{\bfa}_{r(i)}^i\big) -
          \sum_{i\in n} b_i\,x^i\big(r(i),{\bfa}_{r(i)}^i\big) +
          \sum_{i\in n} b_i\,x^i\big(r(i),{\bfa}_{r(i)}^i\big)\Big\|
 \bigg|\ ,\hspace*{5ex}
\end{multline*}
which by the triangle inequality is
\begin{multline*}
\le \sum\limits_{i\in n}|c_i-b_i| \,\big{\|}x^i\big(t(i),
                   {\bfa}_{t(i)}^i\big)\big{\|}
 + \bigg{|}\Big{\|}\sum\limits_{i\in n} b_i\,x^i\big(t(i),
                   {\bfa}_{t(i)}^i\big)\Big{\|} -
                   \Big{\|}\sum\limits_{i\in n} b_i\,x^i\big(r(i),
                   {\bfa}_{r(i)}^i\big)\Big{\|}\bigg{|}\\
\hfill + \sum\limits_{i\in n}|c_i-b_i| \,\big{\|}x^i\big(r(i),
                   {\bfa}_{r(i)}^i\big)\big{\|}\hspace*{1.2ex}\\
 \hspace*{2ex}< n\delta_n + \delta_n + n\delta_n < \eps_n\,,\hfill
\end{multline*}
provided $\delta_n<\frac{\eps_n}{2n+1}$.
\end{proof}

\begin{rmk} One obtains as a limit a norm on $c_{00}$ (the
linear space of finitely supported sequences of scalars),
$\big{\|}\sum_{i\in k}b_i e_i\big{\|}$, where $(e_i)_{i\in\omega}$ is
the unit vector basis for $c_{00}$.
\end{rmk}

We say that $(e_i)_{i\in\omega}$ is a {\bf strong asymptotic model\/}
generated by the strong basic array $(x_i^n)_{n,i\in\omega}$, the
special partition $P\in\langle \omega\rangle^\omega$ and the set of
sequences $\{ {\bfa}_F^i
:i\in\omega$, $F\in [\omega]^{<\omega}\}$. In this case, it is also
easy to see that $(e_i)_{i\in\omega}$ is an asymptotic model of $X$
generated by the basic array $(y_i^n)_{n,i\in\omega}$, where $$y_i^n
= x^n (P(i),{\bfa}^n_{P(i)})\quad \mbox{for}\ n,i\in\omega\,.$$

Thus, asymptotic models can be generated by a stronger type of
convergence. We do not have an application for this. However, it
could prove useful in attacking some of the problems in
Section\;\ref{sec:odd}; those of the type where the assumption is
that every asymptotic model is of a certain type.

We note several special cases of strong asymptotic models
$(e_i)_{i\in\omega}$ generated by $(x_i^n)_{n,i\in\omega}$, $P\in
\spart$ and $\{{\bfa}_F^i: i\in\omega,\;F\in\fin \}$.
\beginsubnr{rmk331}
Let $(x_i)_{i\in\omega}$ be a normalized basic
sequence in $X$ and set $x_i^n=x_i$ for all $n,i\in\omega$. Let
${\bfa}_F^i=(1,0,0,\ldots)$ for all $i\in\omega$ and $F\in\fin$. Then
$(e_i)_{i\in\omega}$ is a spreading model of a subsequence of
$(x_i)_{i\in\omega}$.
\end{subnr}
\beginsubnr{rmk332}
Let $x_i^n=x_i$ for all $n,i\in\omega$, where
again $(x_i)_{i\in\omega}$ is a fixed normalized basic sequence in
$X$. For $i\in\omega$ let ${\bfa}^i$ be a not identically zero
sequence of scalars and set ${\bfa}_F^i={\bfa}^i$ for each
$F\in\fin$. (The non-zero condition is technically violated here, but
we can assume that for some $Q\in\spart$, ${\bfa}_{Q(j)}^i$ is not
identically zero in the first $|Q(j)|$ coordinates if $i\le j$ and
use the theorem to choose $P\in\la Q\ra^\omega$.) In this case we
shall say that $(e_i)_{i\in\omega}$ is a {\bf strong asymptotic model of
$(x_i)_{i\in\omega}$ generated by $P$ and $({\bfa}^i)_{i\in\omega}$}.
\end{subnr}
\beginsubnr{rmk333}
Assume that we are in the situation of
{\small (\ref{rmk332})} with in addition ${\bfa}^i={\bfa}$ for all
$i\in\omega$ and some fixed ${\bfa}$. Then we will say that
$(e_i)_{i\in\omega}$ is an {\bf strong asymptotic model of
$(x_i)_{i\in\omega}$ generated by $P$ and ${\bfa}$}. In this case,
$(e_i)_{i\in\omega}$ is also a spreading model of a normalized block
basis of $(x_i)_{i\in\omega}$.

Indeed, for each $i\in\omega$ let
$y_i=x\big(P(i),{\bfa}\big)$, then $(y_i)_{i\in\omega}$ is a
normalized block basis of $(x_i)_{i\in\omega}$. Also from the
definitions, given $n\in\omega$ and $(b_i)_{i\in n}\in [-1,1]^n$,
$$\bigg{|}
  \Big{\|} \sum_{i\in n} b_i\,y_{j_i} \Big{\|}-
  \Big{\|} \sum_{i\in n} b_i\,e_{i}   \Big{\|}
  \bigg{|} \le \eps_n\,,$$
provided that $n\le j_0<\ldots <j_{n-1}$. Thus, $(e_i)_{i\in\omega}$
is a spreading model of $(y_i)_{i\in\omega}$.
\end{subnr}

\beginsubnr{rmk3-4-3}
If $(e_i)$ is an asymptotic model generated by the strong basic array
$(x_i^n)_{n,i\in\omega}$ then $(e_i)$ is a strong asymptotic model
generated by $(x_i^n)$, $P$ and $({\bfa} _F^i)$ where $P(i) = \{i\}$
and each ${\bfa}_F^i = (1,0,0,\ldots)$.
\end{subnr}

\begin{prp}\label{prp:37}
Let $(e_i)_{i\in\omega}$ be an asymptotic model of $X$ generated by
the basic array $(x_i^n)$. Suppose $(x_i^n)$ is either a weakly null
array or a block basis array (w.r.t.\ some basic sequence in $X$).
Let $(f_i)_{i\in\omega}$ be a normalized block basis of
$(e_i)_{i\in\omega}$. Then $(f_i)_{i\in\omega}$ is also an asymptotic
model of $X$.
\end{prp}

\begin{proof}
Let $(x_i^n)_{n,i\in\omega}$ generate $(e_i)_{i\in\omega}$. Choose
$Q\in\spart$ and ${\bfa}^i$'s such that for every $i\in\omega$,
$|Q(i)|$ is equal to the length of ${\bfa}^i$ and $f_i=e\big(
Q(i),{\bfa}^i \big)$. We shall define a new $K$-basic array
$(y_i^n)_{n,i\in\omega}$ which asymptotically generates
$(f_i)_{i\in\omega}$. For $i\in\omega$ let $\tilde{x}_i$ be the
$i^\te$ diagonal of the array $(x_i^n)_{n,i\in\omega}$, so,
$\tilde{x}_i=(x_i^0,x_{i+1}^1,\ldots,x_{i+n}^n,\ldots)$. As before,
let $\tilde{x}_i\big( F,{\bfa}\big)$ be defined relative to this
sequence. For $n,i\in\omega$ let $z_i^n=\tilde{x}_i\big( Q(n),
{\bfa}^n\big)$. By passing to a subarray of $(z_i^n)_{n,i\in\omega}$
we obtain, as in Lemma\;\ref{lem:normalized}, an array
$(y_i^n)_{n,i\in\omega}$ which is $K$-basic and asymptotically
generates $(f_i)_{i\in\omega}$.
\end{proof}

\begin{rmk}\label{rem:3-5a}
The proposition is false in the general setting. The problem with the
proof is that the rows of $(y_i^n)_{n,i\in\omega}$ need not be
uniformly basic. We sketch how to construct a space $X$ admitting an
asymptotic model $(x_i)_{i\in\omega}$ for which some normalized block
basis $(y_i)_{i\in\omega}$ of $(x_i)_{i\in\omega}$ is not an
asymptotic model of $X$. First we define a norm on
$[(x_i)_{i\in\omega}]$ where $(x_i)_{i\in\omega}$ is a linearly
independent sequence in some linear space. Let $n_i\uparrow \infty$
rapidly and let $(E(i))_{i\in\omega}$ be a special partition of
$\omega$ with $|E(i)|=n_i$. Set for $x= \sum a_i x_i$, $\|x\| = \max
(\|(a_i)\|_{\ell_2}, (\|E_ix\|_{\ell_1})_{T^*})$ where $E_ix$ is the
restriction of $x$ to $E_i$ and $T^*$ is the dual norm to Tsirelson's
space $T$. $(x_i)_{i\in\omega}$ is an unconditional basis for the
reflexive space $[(x_i)_{i\in\omega}]$. Let $y_i = \frac1{|E_i|}
\sum_{j\in E_i}x_j$. Then $(y_i)_{i\in\omega}$ is a normalized block
basis of $(x_i)_{i\in\omega}$ which is equivalent to the unit vector
basis of $T^*$.

Let $X = [(x_i)_{i\in\omega}] \oplus_\infty (\sum\ell_1)_{\ell_2}$.
Let $x_i^n = x_i + e_i^n$ where $(e_i^n)_{i\in\omega}$ is the unit
vector basis of the $n^{\te}$ copy of $\ell_1$ in $(\sum
\ell_1)_{\ell_2}$. Then $(x_i^n)_{n,i\in\omega}$ is a basic array and
generates the asymptotic model $(x_i)_{i\in\omega}$. It can be shown
however that $(y_i)_{i\in\omega}$ is not an asymptotic model of $X$.
The basis $(x_i)_{i\in\omega} \cup (e_i^n)_{n,i\in\omega}$ for $X$ is
boundedly complete and unconditional and thus by passing to a
subarray we may assume that $y_i^n = z_n + w_i^n$ where $z_n\in X$
and $(w_i^n)_{i\in\omega}$ is a seminormalized block basis of the
basis above, in some order, for $X$.

If $P$ is the natural projection of $X$ onto $(\sum
\ell_1)_{\ell_2}$, there must exist $m$ so that, passing to another
subarrary, $\inf_{n\ge m} \inf_{i\ge n} \|P(w_i^n)\|>0$. Otherwise, a
subsequence of $(y_i)_{i\in\omega}$ would be generated by a block
basis array of $(x_i)_{i\in\omega}$ which is impossible. It then
follows that $(y_i)_{i\in\omega}$ must dominate the unit vector basis
of $\ell_2$ due to the structure of $(\sum \ell_1)_{\ell_2}$. Again,
this is false.

It is always true, however, that a normalized block basis of any
spreading model of $X$ is again an asymptotic model of $X$. The
difficulty of choosing $(y_i^n)_{n,i\in\omega}$ to be a basic
subarray (in the proof of Proposition\;\ref{prp:37}) disappears in
this instance.
\end{rmk}

We next collect together a number of remarks and propositions concerning
asymptotic models.

\begin{Obs}\label{Obs:3-6}
\beginsubnr{rmk344}
It is not true in general that an asymptotic model
$(e_i)_{i\in\omega}$ of a basic sequence $(x_i)_{i\in\omega}$ (as in
{\small (\ref{rmk332})}) will be equivalent to a block basis of some
spreading model of $X$, even if $X$ is reflexive.
\jump \hspace*{\parindent} Indeed, consider
$X=\big(\sum\ell_2\big)_{\ell_p}$, with $2<p<\infty$. The only
spreading models of $X$ are $\ell_p$ (isometrically) and $\ell_2$
(isomorphically). This is well-known and easily verified. Letting
$(e_i^n)_{i\in\omega}$ be the unit vector basis of the ``$n^\te$
copy'' of $\ell_2$ in $X$, we can order the unconditional basis
$(e_i^n)_{n,i\in\omega}$ for $X$ as follows:
$$\Big(e_0^0,e_1^0,e_0^1,e_2^0,e_1^1,e_0^2,e_3^0,e_2^1,e_1^2,e_0^3,\ldots
\Big)$$ Take $P(0)=\{0\}$, $P(1)=\{1,2\}$, $P(2)=\{3,4,5\}$,
$P(3)=\{6,7,8,9\}$, $\ldots$ Then this basis along with
$P=\{P(i):i\in\omega\}\in\spart$ generates a strong asymptotic model
$(e_i)_{i\in\omega}$ for the sequence of ${\bfa}^i$'s defined as
follows. Let $n_i$ be positive integers increasing to $\infty$ and
take ${\bfa}^0={\bfa}^1=\ldots ={\bfa}^{n_0}=(1,0,0,0,\ldots)$,
${\bfa}^{n_0+1}=\ldots ={\bfa}^{n_0+n_1}=(0,1,0,0,\ldots)$,
${\bfa}^{n_0+n_1+1}=\ldots ={\bfa}^{n_0+n_1+n_2}=(0,0,1,0,\ldots)$,
{\sl etc}. Then $(e_i)_{i\in\omega}$, as is easily checked, is the
unit vector basis of $\big(\sum\ell_2^{n_i}\big)_{\ell_p}$, which is
not equivalent to a block basis of any spreading model in $X$.
\end{subnr}
\end{Obs}
\beginsubnr{rmk345} One can slightly change the space in {\small
(\ref{rmk344})} to obtain a reflexive space $X$ and a strong asymptotic model
$(e_i)_{i\in\omega}$ which is both not equivalent to a block basis of
a spreading model nor does $E=[(e_i)_{i\in\omega}]$ embed into $X$.
The same sort of scheme as presented in {\small (\ref{rmk344})} works
for $X=\big(\sum T\big)_{\ell_2}$, the $\ell_2$ sum of Tsirelson's
space $T$ (see\;\cite{FigielJohnson}). The only spreading models of
this space are all isomorphic to $\ell_1$ or $\ell_2$.
For, if $P_n$ is the norm $1$ natural projection of $X$ onto the ``$n^\te$
copy'' of $T$ in $X$, and $(x_i)_{i\in\omega}$ is a normalized basic
sequence in this reflexive space, then passing to a subsequence we
may assume {\sl either\/}: for all $n$, $\lim_{i\to\infty}\| P_n
x_i\|=0$, in which case, by a gliding hump argument,
$(x_i)_{i\in\omega}$ has $\ell_2$ as a spreading model; {\sl or\/}:
for some $n$, $\lim_{i\to\infty}\| P_n x_i\|>0$, in which case
$(x_i)_{i\in\omega}$ has a subsequence whose spreading model is
isomorphic to $\ell_1$. Now, if we use the basis ordering of {\small
(\ref{rmk344})} and the same $P(i)$'s, and take the ${\bfa}^i$'s to
be such that for each sequence $(0,0,\ldots,0,1,0,0,\ldots)$,
infinitely many ${\bfa}^i$'s are equal to this sequence, then we
obtain $\big(\sum\ell_1 \big)_{\ell_2}$ as a strong asymptotic model. This
does not embed into $X$.
\end{subnr}
\beginsubnr{rmk346} Spreading models join the infinite and
arbitrarily spread out and finite dimensional structure of $X$.
Another such joining is the theory of {\em asymptotic structure\/}
(see\;\cite{Maurey.etal}). In its simplest form this can be described
as follows. Suppose $X$ has a basis $(x_i)_{i\in\omega}$. For a
positive $n\in\omega$, a normalized basic sequence $(e_i)_{i\in n}$
belongs to the $n^\te$-asymptotic structure of $X$, denoted
$\{X\}_n$, if for all $\eps>0$, given $m_0\in\omega$ there exists
$y_0\in S_{\la (x_i)_{i\in\omega\setminus m_0}\ra}$, so that for all
$m_1\in\omega$ there exists $y_1\in S_{\la (x_i)_{i\in\omega\setminus
m_1}\ra}$,\,$\ldots$, so that for all $m_{n-1}\in\omega$ there exists
$y_{n-1}\in S_{\la (x_i)_{i\in\omega\setminus m_{n-1}}\ra}$, so that
$(y_i)_{i\in n}$ is $(1+\eps)$-equivalent to $(e_i)_{i\in n}$. (Here,
$S_{\la (x_i)_{i\in\omega\setminus m_{j}}\ra}$ denotes the unit
sphere of the linear span of $\{x_i:i\in\omega\setminus m_j\}$.)\jump
One difference between this and spreading models is that spreading
models are infinite. However one can paste together the elements of
the sets $\{X\}_n$ as follows.  $(e_i)_{i=1}^\infty$ is an {\bf
asymptotic version\/} of $X$ if for all $n$, $(e_i)_{i=1}^n \in
\{X\}_n$ \cite{Maurey.etal}. But certain infinite threads are lost
nonetheless. Furthermore, spreading models arise from ``every
normalized basic sequence has a subsequence$\ldots$''. $\{X\}_n$ can
be described in terms of infinitely branching trees of length $n$.
The initial nodes and the successors of any node form a normalized
block basis of $(x_i)_{i\in\omega}$. We can label such a tree as
$T_n=\{x_{(m_0,\ldots,m_k)}: 0\le m_0<\ldots <m_k,\;k\in n\}$ ordered
by $x_\alpha\le x_\beta$ if the sequence $\alpha$ is an initial
segment of $\beta$. Then $(e_i)_{i\in n}\in\{X_n\}$ iff there exists
a tree $T_n$ so that for all $\eps >0$ there exists $n_0$ so that if
$n_0\le m_0<\ldots <m_{n-1}$, then $\big(x_{(m_0,\ldots,
m_k)}\big)_{k\in n}$ is $1+\eps$-equivalent to $(e_i)_{i\in n}$. This
stronger structure yields in some sense a more complete theory than
that of spreading models where a number of problems remain open. The
theory of asymptotic models generated by block basic arrays, while
being closer to that of spreading models, lies somewhere between the
two. The theory and open problems of spreading models and asymptotic
structure motivate some of our questions and results below.\jump
Further, it is clear that if $X$ has a basis $(x_i)_{i\in\omega}$ and
$(e_i)_{i\in\omega}$ is an asymptotic model of $X$ generated by a
block basis array (w.r.t.\ $(x_i)_{i\in\omega}$), then for all $n$,
$(e_i)_{i\in n}\in\{X\}_n$.
\end{subnr}
\beginsubnr{rmk347} Suppose that $X$ has a basis and that all spreading
models of a normalized block basis are equivalent. {\sl Must all
spreading models be equivalent to the unit vector basis of $c_0$ or
$\ell_p$ for some $1\le p<\infty$\,?\/} This question, due to Spiros
Argyros, remains open. Some partial results are in
\cite{Androulakis.etal}. The analogous question for asymptotic models
has a positive answer.\jump \hspace*{\parindent} Indeed, suppose that
all asymptotic models of all block basis arrays of $X$ are
equivalent. If $(\tilde{e}_i)_{i\in\omega}$ is a spreading model of
such a space, then all of its normalized block bases, being
asymptotic models by Proposition\;\ref{prp:37}, must be equivalent
and the result follows from Zippin's Theorem (see \cite{Zippin} or
\cite[p.\,59]{LindenstraussTzafriri}).
\end{subnr}
\beginsubnr{rmk348} If $X$ is reflexive and $(e_i)_{i\in\omega}$ is an
asymptotic model of $X$, then $(e_i)_{i\in\omega}$ is suppression-$1$
unconditional. More generally, this holds if $(e_i)_{i\in\omega}$ is
generated by $(x_i^n)_{n,i\in\omega}$ where for each $n\in\omega$,
$(x_i^n)_{i\in\omega}$ is weakly null. \jump The proof is very much
the same as the analogous result for spreading models. Let
$(b_i)_{i\in n}\in [-1,1]^n$ and $i_0\in n$. We need only show
$\big{\|}\sum_{i\in n\setminus \{i_0\}} b_i e_i\big{\|}\le
\big{\|}\sum_{i\in n} b_i e_i\big{\|}$.
\end{subnr}

Let $m\ge n$. Since $(x_j^{i_0})_{j\in\omega}$ is weakly null there
exists a convex combination  of small norm: $\|\sum_{p\in k} c_p
x_{m+i_0+p}^{i_0}\| <\eps_m$. For $p\in k$ we consider the vector
$$y_p = \sum_{i\in i_0} b_i x_{m+i}^i + b_{i_0} x_{m+i_0+p}^{i_0} +
\sum_{i= i_0+1}^{n-1} b_i x_{m+k+i}^i\,.$$ $|\,\|\sum_{i\in n}  b_i
e_i\| - \|y_p\|\,| <\eps_m$ and so $$\bigg\| \sum_{p\in k} c_p
y_p\bigg\| \le \bigg\| \sum_{i\in n} b_i e_i\bigg\| + \eps_m\,.$$ but
also $$\bigg\| \sum_{p\in k} c_p y_p\bigg\| \ge \bigg\|
\sum_{\substack{i\in n\\ i\ne i_0}} b_i e_i\bigg\| - \eps_m -
|b_{i_0}| \eps_m$$ and this yields the desired inequality.

\beginsubnr{rmk349} In general, the $n^\te$ asymptotic structure
$\{X\}_n$ of a Banach space $X$ with a basis $(x_i)_{i\in\omega}$
does not coincide with $\big{\{} (e_i)_{i\in n}: (e_i)_{i\in\omega}$
is an asymptotic model generated by a block basis array of
$(x_i)_{i\in\omega}\big{\}}$.
In fact, these may be vastly different
for every subspace of $X$ generated by a block basis of
$(x_i)_{i\in\omega}$.\jump To see this we recall that in
\cite[Section\;3]{OdellSchlumprecht1}, a reflexive $X$ is constructed
so that $(y_i)_{i\in n}\in\{X\}_n$ for all normalized monotone basic
sequences $(y_i)_{i\in n}$. Since this includes the highly
unconditional summing basis (of length $n$) the claim follows from
{\small (\ref{rmk348})}.
\end{subnr}

\beginsubnr{rmk381} It is possible for a space $X$ to have $\ell_1$ as
an asymptotic model yet no spreading model of $X$ is isomorphic to
$\ell_1$, nor to $c_0$ or any $\ell_p$ ($1<p<\infty$).\jump Indeed,
the reflexive
space $X$ constructed in \cite{Androulakis.etal} has the property
that no spreading model is isomorphic to $\ell_p$ ($1\le p < \infty$)
nor $c_0$. Yet every spreading model of $X$ contains an isomorphic
copy of $\ell_1$. \end{subnr}
\beginsubnr{rmk382} There exists a reflexive space $X$ for which no asymptotic
model contains an isomorphic copy of $c_0$ or $\ell_p$ ($1\le
p\le\infty$). \jump $X$ is the space constructed in
\cite{OdellSchlumprecht95}; we recall the example: $\|\cdot\|$ is a
norm on $c_{00}$ satisfying the following implicit equation. $$\| x\|
:=\max\bigg{\{} \|x\|_{c_0},\Big( \sum_{k\in\omega}
\|x\|_{n_k}^2\Big)^{1/2}\bigg{\}}\,,$$ where $\|x\|_{n_k}
=\sup\Big{\{}\frac{1}{f(n_k)} \sum_{i\in n_k}\| E_i x\|: E_0<\ldots <
E_{n_k-1}\Big{\}}$,\; $f(n_k)=\log_2(1+n_k)$ and $(n_k)_{k\in\omega}$
is a sequence of positive integers satisfying
$\sum_{k\in\omega}\frac{1}{f(n_k)}<\frac{1}{10}$. $X$ is the
completion of $c_{00}$ under this norm. The unit vector basis
$(u_i)_{i\in\omega}$ of $c_{00}$ is a $1$-unconditional basis for $X$
and $X$ is reflexive. The fact that $X$ does not admit an asymptotic
model $(e_i)_{i\in\omega}$ equivalent to the unit vector basis of
$\ell_1$ (and hence, by Proposition\;\ref{prp:37}, no asymptotic
model $E$ contains $\ell_1$) is similar to the proof in
\cite{OdellSchlumprecht95} that no spreading model is isomorphic to
$\ell_1$, and so we shall only sketch the argument.\jump Suppose that
$(e_i)_{i\in\omega}$ is an asymptotic model of $X$ and is equivalent to
the unit vector basis of $\ell_1$.
We may assume that $(e_i)_{i\in\omega}$ is generated by
the basic array $(x_i^n)_{n,i\in\omega}$ where each
$(x_i^n)_{i\in\omega}$ is a normalized block basis of
$(u_i)_{i\in\omega}$. By iteratively passing to a subsequence of each
row $(x_i^n)_{i\in\omega}$ and diagonalizing, we may assume that
$\big( {\|} x_j^n {\|}_{n_i} \big)_{i\in\omega}$ converges weakly in
$B_{\ell_2}$ as $j\to\infty$ to ${\bfa}^n\in B_{\ell_2}$. Considering
the sequence $({\bfa}^n)_{n\in\omega}\subs B_{\ell_2}$ and passing to
a subsequence of the rows, we may assume that
$({\bfa}^n)_{n\in\omega}$ converges weakly in $B_{\ell_2}$ to some
${\bfa}\in B_{\ell_2}$. This corresponds to passing to a subsequence
of $(e_i)_{i\in\omega}$, but that is still equivalent to the unit
vector basis of $\ell_1$ and so we lose nothing here. Thus, we are in
the situation where the limit distribution in $\ell_2$ of the $n^\te$
row $(x_i^n)_{i\in\omega}$ is ${\bfa}^n$ and therefore we can assume
$\big( \| x_i^n {\|}_{n_j} \big)_{j\in\omega}$ in $\ell_2$ is equal
to ${\bfa}^n + {\boldsymbol{h}}_i^n$, where $({\boldsymbol{h}}_i^n
)_{i\in\omega}$ is weakly null in $\ell_2$. Furthermore, ${\bfa}^n
={\bfa}+ {\boldsymbol{h}}^n$, where ${\boldsymbol{h}}^n$ is weakly
null in $\ell_2$ and hence, we may assume, a block basis in $\ell_2$.
In this manner, for any $N$ and $(b_i)_{i\in N}\in [-1,1]^N$ we have
$\big{\|} \sum_{i\in N} b_i e_i \big{\|} \approx \big{\|} \sum_{i\in
N} b_i x_{k_i}^i \big{\|}$, provided $N\le k_0<\ldots
<k_{N-1}$.\hfill\\
Now we can also assume that $\big{\|} \sum_{i\in
N} b_i e_i \big{\|} \ge 0.99\cdot\sum_{i\in N} |b_i|$. This is
because $\ell_1$ is not distortable (see \cite{James}) and every
block basis of an asymptotic model of $X$ is (by Proposition\;\ref{prp:37})
also an asymptotic model. Thus, by carefully choosing the $k_i$'s, we
have $ 0.99\cdot\sum_{i\in N} |b_i| < \big{\|} \sum_{i\in N} b_i
x_{k_i}^i \big{\|}$ where $\big( \| x_{k_i}^i \| \big)_{\ell_2}
\approx {\bfa}^i + {\boldsymbol{h}}^i + {\boldsymbol{h}}_{k_i}^i$ and
the vectors $\big( {\boldsymbol{h}}^i + {\boldsymbol{h}}_{k_i}^i
\big)_{i\in N}$ are a block basis in $\ell_2$. At this point, we use
the argument in Theorem\;1.3 of \cite{OdellSchlumprecht95} to see
that, if $N$ is sufficiently large depending upon ${\bfa}$, this is
impossible.\hfill\\ Furthermore, the arguments of
\cite{OdellSchlumprecht95} apply easily to show that it is not
possible to have an asymptotic model $(e_i)_{i\in\omega}$ equivalent
to the unit vector basis of $c_0$ or $\ell_p$ ($1<p<\infty$), which
completes the proof of {\small (\ref{rmk382})}.
\end{subnr}
\beginsubnr{rmk383} The proof of {\small (\ref{rmk382})} actually reveals that
no spreading model of an asymptotic model of $X$ can be isomorphic to
$\ell_1$ (or $c_0$ or any $\ell_p$). For if
$E=\big{[}(e_i)_{i\in\omega} \big{]}$ is an asymptotic model of $X$,
then any spreading model of $E$ is necessarily a spreading model of a
normalized block basis $(f_i)_{i\in\omega}$ of $(e_i)_{i\in\omega}$
and this in itself is an asymptotic model of $X$. Let
$(\tilde{e}_i)_{i\in\omega}$ be the spreading model of
$(f_i)_{i\in\omega}$. The proof shows that, for sufficiently large
$N$, we cannot have $\big{\|} \sum_{i\in N} b_i f_{k_i} \big{\|} \ge
0.99\cdot\sum_{i\in N} |b_i|$ for all $(b_i)_{i\in N}\in [-1,1]^N$
and any $k_0<\ldots
<k_{N-1}$.
\end{subnr}
\beginsubnr{rmk384}
In \cite{Androulakis.etal} a reflexive Banach space $X$ is
constructed for which no spreading model is reflexive, isomorphic to
$c_0$ or isomorphic to $\ell_1$. However, every $X$ admits an
asymptotic model which is either reflexive or isomorphic to $c_0$ or
$\ell_1$.\jump Indeed, $X$ admits a spreading model $\tilde{E}$ with
an unconditional basis and by \cite{James}, $\tilde{E}$ is either
reflexive or contains an isomorphic copy of $c_0$ or $\ell_1$. So,
the result follows by Remark\;\ref{rem:3-5a}.
\end{subnr}

There is a big difference between considering all asymptotic models
of $X$ and of those generated by weakly null basic arrays or block
basic arrays as our next proposition illustrates. Also it illustrates
again the difference between the class of spreading models and
asymptotic models:  if $(e_i)_{i\in\omega}$ is a spreading model of
$c_0$, then $(e_i)_{i\in\omega}$ is equivalent to either the summing
basis or the unit vector basis of $c_0$.

\begin{prp}\label{prp:3-7}
Let $(e_i)_{i\in\omega}$ be a normalized bimonotone basic sequence.
Then $(e_i)_{i\in\omega}$ is $1$-equivalent to an asymptotic model of
$c_0$.
\end{prp}

\begin{proof}
Let $\eps_k\downarrow 0$. For all positive integers $k$ there exist
$n_k\in\omega$ and vectors $(x_i^k)_{i\in k} \in
S_{\ell_\infty^{n_k}}$ so that $$(1-\eps_k) \bigg\| \sum_{i\in k} a_i
e_i \bigg\| \le \bigg\| \sum_{i\in k} a_i x_i^k\bigg\| \le \bigg\|
\sum_{i\in k} a_i e_i\bigg\|$$ for all $(a_i)_{i\in k}\in {\mathbb
R}^k$. Indeed, we choose $(f_i^k)_{i\in n_k} \subseteq
B_{[(e_i)_{i\in\omega}]^*}$ so that $\sup_{i\in n_k} |f_i^k (e) | \ge
(1-\eps_k)\|e\|$ for $e\in [(e_i)_{i\in k}]$ and $f_i^k(e_i) =1$ for
$i\in k$, and let $x_i^k = T^k e_i$, where $T^k :[(e_i)_{i\in k}]\to
\ell_\infty^{n_k}$ is given by $T^k e= \big(f_i^k(e)\big)_{i\in
n_k}$.

We write $c_0 = \big(\sum \ell_\infty^{n_k}\big)_{c_0}$ and regard
$(x_i^k)_{i\in k}$ as being contained in the indicated copy of
$\ell_\infty^{n_k} \subseteq c_0$. Let $(y_i^k)_{k\in\omega,\, i\ge
k}$ be defined by $y_i^0 =  x_0^1 +\cdots + x_0^{i+1}$ and in general
$y_i^k =  x_k^{k+1} +\cdots + x_k^{i+1}$.

It is easy to check that $(y_i^k)_{k,i\in\omega}$ is a basic array
(the rows are equivalent to the summing basis) and this array
generates $(e_i)_{i\in\omega}$.
\end{proof}

\begin{rmk}\label{rem:DLT}
Recall \cite{DLT00} that a basic sequence $(x_i)_{i\in\omega}$ is
said to be {\it asymptotically isometric to ${c_0}$}, if for some
sequence $\eps_n \downarrow 0$ for all $(a_n)_{n\in\omega}\in c_0$,
$$\sup_n (1-\eps_n) |a_n| \le \Big\| \sum_{n\in\omega} a_n x_n\Big\|
\le \sup_n |a_n|\ .$$ In this case the proof of
Proposition\;\ref{prp:3-7} can be adopted to yield that
$[(x_i)_{i\in\omega}]$ admits all normalized bimonotone basic
sequences as asymptotic models. In general, using that $c_0$ is not
distorable \cite{James}, one has that if $X$ is isomorphic to $c_0$
then for all $K>1$ there exists $C(K)$ so that if
$(e_i)_{i\in\omega}$ is a normalized $K$-basic sequence, then $X$
admits an asymptotic model $C(K)$-equivalent to $(e_i)_{i\in\omega}$.
We do not know if the conclusion to Proposition\;\ref{prp:3-7} holds
in this case. We also do not know if this property characterizes
spaces containing $c_0$ (see the open problems in
Section\;\ref{sec:odd}). By way of contrast it is easy to see that
all asymptotic models of $\ell_p$ $(1<p<\infty)$ are 1-equivalent to
the unit vector basis of $\ell_p$. Moreover we have
\end{rmk}

\begin{prp}\label{prp:3-8}
If $(e_i)_{i\in\omega}$ is an asymptotic model of $\ell_1$ then
$(e_i)_{i\in\omega}$ is equivalent to the unit vector basis of
$\ell_1$.
\end{prp}

\begin{proof}
Let $(x_i^n)_{n,i\in\omega}$ be a $K$-basic array generating
$(e_i)_{i\in\omega}$. Since each row is $K$-basic there exists
$\delta>0$ so that for all $n,m\in\omega$ there exists $k\in\omega$ with
$\|P^m (x_i^n)\|>\delta$ for $i\ge k$ where $P^m$ is the tail
projection of $\ell_1$, $P^m (a_i) = (0,\ldots,0,a_m,a_{m+1},\ldots)$.
Using that the unit vector basis
of $\ell_1$ is boundedly complete we can find a subsequence
$(y_i^n)_{i\in\omega}$ of each row $(x_i^n)_{i\in\omega}$ of the form
$y_i^n = y_n +h_i^n$ where $h_i^n\to 0$ weak* in $\ell_1$ as
$i\to\infty$ and $\|h_i^n\| \ge\delta$. Thus, up to arbitrarily small
perturbations, we may assume $h_i^n$ and $h_j^n$ are disjointly
supported for $i\ne j$. And doing all this by a diagonal process we
can assume that $(y_i^n)_{n,i\in\omega}$ is a subarray of
$(x_i^n)_{n,i\in\omega}$. It follows easily that $$\big{\|} \sum a_i
e_i\big{\|} \ge \delta \sum |a_i|\,.$$\vspace*{-4ex}

\end{proof}
\smallskip

{From} Proposition\;\ref{prp:3-7} we see that $\ell_1$ can be an
asymptotic model of a space $X$ with a basis without being an
asymptotic model generated by a block basic array. But this cannot
happen in a boundedly complete situation:

\begin{prp}\label{prp:3-9}
Let $(d_i)_{i\in\omega}$ be a boundedly complete basis for $Y$ and
let $X\subseteq Y$ be a weak* closed subspace. If $\ell_1$ is an
asymptotic model of $X$, then $\ell_1$ is an asymptotic model
generated by a basic array $(x_i^n)_{n,i\in\omega}$ where for each
$n$, $(x_i^n)_{i\in\omega}$ is weak* null.
\end{prp}

In this proposition the weak* topology on $Y$ is the natural one
generated by regarding $Y$ as the dual space of
$[(d_i^*)_{i\in\omega}]$, where the $d_i^*$'s are the biorthogonal
functionals of the $d_i$'s (this is, for all $i,j$, $d_i^* d_j=
\delta_j^i$). Thus, $d_n = \sum a_i^n d_i \to d = \sum a_i d_i$ weak*
if $(d_n)_{n\in\omega}$ is bounded and $a_i^n\to a_i$ for each
$n\in\omega$.

\begin{proof}[Proof of Proposition\;\ref{prp:3-9}]
Let $(y_i^n)_{n,i\in\omega}\subseteq X$ generate the asymptotic model
$(e_i)_{i\in\omega}$ which is equivalent to the unit vector basis of
$\ell_1$. As in the preceding proposition by passing to a subarray we
may assume $y_i^n = f^n + x_i^n$ where for each $n$,
$(x_i^n)_{i\in\omega}$ is weak* null and $(f^n)_{n\in\omega}\subseteq
X$. If $(f^n)_{n\in\omega\setminus k}$ is not equivalent to the unit
vector basis of $\ell_1$ for some $k$, then some block sequence of
absolute convex combinations of the $f^n$'s is norm null. We use this
(as in the proof of Proposition\;\ref{prp:37}) to generate a new
basic array of the same form where $\|f^n\| < \eps_n$ for
$\eps_n\downarrow 0$ rapidly, and so, a subarray of
$\Big(\frac{x_i^n}{\|x_i^n\|}\Big)_{n,i\in\omega}$ generates the unit
vector basis of $\ell_1$.
\end{proof}

The asymptotic models of $L_p$ $(1<p<\infty)$ are necessarily unconditional
and in fact every normalized unconditional basic sequence in $L_p$ is
equivalent to an asymptotic model.

\begin{prp}\label{newprop3.12}
Let $1<p<\infty$. There exists $K_p <\infty$ so that if
$(x_i)_{i\in\omega}$ is a normalized $K$-unconditional basic sequence
in $L_p$ then $(x_i)_{i\in\omega}$ is $KK_p$-equivalent to some
asymptotic model of $L_p$.
\end{prp}

\begin{proof}
This follows easily from arguments of Gideon Schechtman \cite{S74}.
There exists $K_p <\infty$ so that $(x_i)_{i\in\omega}$ is
$KK_p$-equivalent to a normalized block basis $(y_i)_{i\in\omega}$ of
the Haar basis $(h_i)_{i\in\omega}$ for $L_p$. Furthermore if
$(z_i)_{i\in\omega}$ is a block basis of $(h_i)_{i\in\omega}$ with
$|z_i| =|y_i|$ for all $i$, then $(x_i)_{i\in\omega}$ is
$KK_p$-equivalent to $(z_i)_{i\in\omega}$. For $n\in\omega$, let
$(y_i^n)_{i\in\omega}$ be a normalized block basis of
$(h_i)_{i\in\omega}$ with $|y_i^n| = |y_n|$ for all $i$. By
Lemma~\ref{lem:normalized}, some subarray of $(y_i^n)_{n,i\in\omega}$
is thus a block basis array of $(h_i)_{i\in\omega}$. By our above
remarks and Propostion~\ref{prp:K-basic}, some subarray of
$(y_i^n)_{n,i\in\omega}$ generates an asymptotic model
$KK_p$-equivalent to $(x_i)_{i\in\omega}$.
\end{proof}

Another natural question is if $X$ has $Y$ as an asymptotic model and
$Y$ has $Z$ as an asymptotic model, does $X$ have an asymptotic model
isomorphic to $Z$\,? If one replaces ``asymptotic model'' in the
question with ``spreading model'', the answer is negative (see
\cite{BeauzamyMaurey}). In the following, we present an example that
shows the answer also to be negative in a strong way for asymptotic models.

\begin{exa}\label{exa} There exist reflexive Banach spaces $X$ and $Y$ so
that $Y$ is a spreading model of $X$, $\ell_1$ is a spreading model
of $Y$ and $\ell_1$ is not isomorphic to any asymptotic model of $X$.
\end{exa}

\begin{proof}
 $X$ and $Y$ will both be completions of $c_{00}$ under certain
norms which will make the unit vector basis of $c_{00}$ an
unconditional basis for each space. We will denote these bases by
$(v_i)_{i\in\omega}$ for $X$ and $(u_i)_{i\in\omega}$ for $Y$. Both
spaces will be reflexive.\jump\hspace*{\parindent} First we construct
the spaces $Y$ and $X$. The construction bears some similarity with
those in \cite{MaureyRosenthal} and
\cite[p.\,123]{LindenstraussTzafriri}. To begin, let
$(m_j)_{j\in\omega}$ be an increasing sequence of integers with
$m_0=1$ and for any $k\in\omega$: $m_0+\ldots +m_k < 2m_k$,
$\sum_{n\in\omega\setminus \{0\}}  \frac1{\sqrt{m_n}} <1$ and
$\frac{(2m_k)^2}{\sqrt{m_{k+1}}} <1$. Let $\cF$ be the subset of
$c_{00}$ given as follows:
\begin{multline*}\cF:=\bigg{\{} f=\sum_{j\in
n}\frac{1_{E_{i_j}}}{\sqrt{m\indd}_{\mindd i_j}} :
n\in\omega,\;|E_{i_j}|\le m_{i_j},\;n\le i_0<\ldots <i_{n-1}\
\text{and}\\ E_{i_k}\cap E_{i_l} =\emptyset\ \text{whenever $k\neq
l$}\,\bigg{\}}\,,\end{multline*}
where $1_{E_{i_j}}\in c_{00}$ is the indicator function,
$$1_{E_{i_j}}(k)=\begin{cases} 1 & \text{if $k\in
E_{i_j}$,}\\
                               0 & \text{otherwise.}
                 \end{cases}$$
For $x\in c_{00}$, let
\begin{multline*}\| x\|_Y:= \sup \bigg{\{}
\Big( \sum_{k\in m} \la f_k,x\ra^3\Big)^{\frac13}:
m\in\omega,\;(f_k)_{k\in m}\subs\cF\ \text{and}\\ \text{the $f_k$'s
are disjointly supported}\,\bigg{\}}\,,
\end{multline*}
where $\la f_k,x\ra$ is the scalar product of $f_k$ and $x$. We say
$E\in\fin$ is {\bf admissible\/} if $\min (E)\ge |E|$ and $g\in
c_{00}$ is admissible if $\supp (g)$ (the support of $g$) is
admissible. Set $\cG:=\{f|_E: E\ \text{is admissible and
$f\in\cF$}\}=\{f\in\cF: f\ \text{is admissible}\}$, and for $x\in
c_{00}$, let
\begin{multline*}\| x\|_X:= \sup \bigg{\{}
\Big( \sum_{k\in m} \la g_k,x\ra^3\Big)^{\frac13}:
m\in\omega,\;(g_k)_{k\in m}\subs\cG\ \text{and}\\ \text{the $g_k$'s
are disjointly supported}\,\bigg{\}}\,.
\end{multline*}
We will also write $g(x)$ for $\langle g,x\rangle$. It is clear that
$(v_j)_{j\in\omega}$ and $(u_j)_{j\in\omega}$ are each
\hbox{suppression-$1$} unconditional bases for $X$ and $Y$,
respectively. Because each basis admits a lower $\ell_3$ estimate on
disjointly supported vectors, neither space contains
$\ell_\infty^n$'s uniformly (see \cite{Johnson}). Thus, both bases
are boundedly complete. Also both bases are shrinking and hence, $X$
and $Y$ are reflexive. To see this for $Y$ (the proof for $X$ is
similar) suppose $(y_i)_{i\in\omega}$ is a normalized block basis of
$(u_j)_{j\in\omega}$ which is not weakly null. By the definition of
the norm in $Y$, and passing to a subsequence of
$(y_i)_{i\in\omega}$, we obtain $f\in\cF$ and $\eps>0$ with $|\la
f,y_j\ra |>\eps$ for all $j$, which is clearly impossible.
\jump\hspace*{\parindent} The sequence $(u_j)_{j\in\omega}$ is
\hbox{$1$-symmetric} and is the spreading model of
$(v_j)_{j\in\omega}$ (since if one moves a vector far enough to the
right in $c_{00}$, then the $Y$ norm expressions all become
allowable). \jump\hspace*{\parindent} Let $E_0<\ldots <E_j<\ldots$ be
sets of natural numbers with $|E_j|=m_j$ and let
$y_j=\frac{1_{E_j}}{\sqrt{m\ind}_{\mind j}}$ (for $j\in\omega$). Then
$\| y_j \|_Y\ge 1$ and $\sup_{j\in\omega} \| y_j\|_Y<\infty$. Indeed,
for some fixed $q\in\omega$, let
$y=\frac{1_{E_q}}{\sqrt{m\ind}_{\mind q}}$. First suppose $f\in\cF$,
and therefore, $f$ is of the form $f=\sum\limits_{j\in n}
\frac{1_{E_{i_j}}}{\sqrt{m\indd}_{\mindd i_j}}$ (for some disjoint
collection $(E_{i_j}) \subseteq \fin$ with $|E_{i_j}|\le m_{i_j}$ and
$n\le i_0 < \cdots < i_{n-1}$). We shall estimate $\la f,y\ra$ from
above, and thus we may assume $\supp (f)\subs E_q$. Write $f=f^1+f^2
+f^3$, where

\begin{align*} f^1=&\sum\limits_{\begin{subarray}{c}
                   j\in n\\
                   i_j<q
                   \end{subarray}}
                   \frac{1_{E_{i_j}}}{\sqrt{m\indd}_{\mindd
                   i_j}}\,,
                   & f^2=&\begin{cases}
      \displaystyle{\frac{1_{E_{q}}}{\sqrt{m\ind}_{\mind q}}} & \text{if some $i_j=q$},\\
      0 & \text{otherwise},
      \end{cases}
                   &
                   f^3=&\sum\limits_{\begin{subarray}{c}
                   j\in n\\
                   i_j>q
                   \end{subarray}}
                   \frac{1_{E_{i_j}}}{\sqrt{m\indd}_{\mindd
                   i_j}}\,.
                   \end{align*}
\noindent By the properties of the sequence $(m_j)_{j\in\omega}$ we
have
 \begin{align*}\la f^1,y\ra & =
                   \displaystyle{\sum\limits_{\begin{subarray}{c}
                              j\in n\\
                              i_j<q
                        \end{subarray}}
                   \frac{|E_{i_j}|}{\sqrt{m\indd}_{\mindd i_j}
                   \sqrt{m\ind}_{\mind q}}\le
                   \frac{2m_{q-1}}{\sqrt{m\ind}_{\mind q}}}\,, &
                \la f^2,y \ra & \le
                   \displaystyle{\frac{m_q}{\sqrt{m\ind}_{\mind q}
                           \sqrt{m\ind}_{\mind q}}=1}\\
\intertext{and}
                \la f^3,y\ra & =
                   \displaystyle{\sum\limits_{\begin{subarray}{c}
                              j\in n\\
                              i_j>q
                        \end{subarray}}
                   \frac{|E_{i_j}|}{\sqrt{m\indd}_{\mindd i_j}
                               \sqrt{m\ind}_{\mind q}}\le
                   \frac{\sqrt{m_q}}{\sqrt{m\inddd}_{\minddd q+1}}}\,.
\end{align*}
Now suppose that $f_k=\sum_{j\in
n_k}\frac{1_{E_{i_j^k}}}{\sqrt{m\indd}_{\mindd i_j^k}}\in \cF$ and
the $(f_k)_{k\in m}$ are disjointly supported with $\supp (f_k)\subs
E_q$ for each $k\in m$. As above, each $f_k$ is of the form
$f_k=f_k^1+ f_k^2+ f_k^3$. Thus by the triangle inequality in
$\ell_3$, $$\Big( \sum_{k\in m}\la f_k,y\ra^3\Big)^{\frac13}\le \Big(
\sum_{k\in m}\la f_k^1,y\ra^3\Big)^{\frac13} + \Big( \sum_{k\in m}\la
f_k^2,y\ra^3\Big)^{\frac13} + \Big( \sum_{k\in m}\la
f_k^3,y\ra^3\Big)^{\frac13}\,.$$ The first term is $\bigg(\sum_{k\in
m} \sum_{\begin{subarray}{c}
                                           j\in n_k\\
                                           i_j^k<q
                                          \end{subarray}}
  \Big(\frac{|E_{i_j^k}|}{\sqrt{m\indd}_{\mindd i_j^k}\sqrt{m\ind}_{\mind q}}
  \Big)^3 \bigg)^{1/3}$,
and since (by the earlier calculation)
 $\sum_{\begin{subarray}{c}
                              j\in n_k\\
                              i_j^k<q
               \end{subarray}}
  \frac{|E_{i_j^k}|}{\sqrt{m\indd}_{\mindd i_j^k}
                     \sqrt{m\ind}_{\mind q}}
  \le\frac{2m_{q-1}}{\sqrt{m_q}}$,
the first term is
  $$
  \le  \displaystyle{\bigg( \Big( \frac{2m_{q-1}}{\sqrt{m_q}}\Big)^2
  \Big( \sum_{k\in m} \sum\limits_{\begin{subarray}{c}
                                          j\in n_k\\
                                          i_j^k<q
                                   \end{subarray}}
 \frac{|E_{i_j^k}|}{\sqrt{m\indd}_{\mindd i_j^k}\sqrt{m\ind}_{\mind q}} \Big)
 \bigg)^{\frac13}}
 \le  \displaystyle{\bigg( \Big( \frac{2m_{q-1}}{\sqrt{m\ind}_{\mind q}}\Big)^2
 \frac{m_q}{\sqrt{m\ind}_{\mind q}}\bigg)^{\frac13} = \bigg( \frac{(2m_{q-1})^2}
 {\sqrt{m\ind}_{\mind q}} \bigg)^{\frac13} <1\,.}
 $$
The second term is of the form
$(\sum_{k\in m} (\frac{l_k}{\sqrt{m_q}\sqrt{m_q}})^3)^{1/3}$,
where $\sum_{k\in m}l_k\le m_q$, and therefore,
it is $\le \sum_{k\in m}\frac{l_k}{m_q}\le 1$.
\jump
The third term is
  $$
  \bigg( \sum_{k\in m} \Big( \sum\limits_{\begin{subarray}{c}
                              j\in n_k\\
                              i_j^k>q
               \end{subarray}}
  \frac{|E_{i_j^k}|}{\sqrt{m\indd}_{\mindd i_j^k}
  \sqrt{m\ind}_{\mind q}}
  \Big)^3 \bigg)^{\frac13}
  \le \sum_{k\in m} \sum\limits_{\begin{subarray}{c}
                              j\in n_k\\
                              i_j^k>q
               \end{subarray}}
  \frac{|E_{i_j^k}|}{\sqrt{m\indd}_{\mindd i_j^k} \sqrt{m\ind}_{\mind
  q}}
  \le \frac{m_q}{\sqrt{m\inddd}_{\minddd q+1}\sqrt{m\ind}_{\mind q}}
  = \frac{\sqrt{m\ind}_{\mind q}}{\sqrt{m\inddd}_{\minddd q+1}} <
  1\,.
  $$

Thus, $(y_j)_{j\in\omega}$ is a seminormalized block basis of
$(u_j)_{j\in\omega}$ in $Y$.
Moreover, from the definition of the
norm, namely $\cF$, if $n\le i_0<\ldots <i_{n-1}$ and $(b_i)_{i\in
n}$ are scalars, then
$\| \sum_{j\in n} b_i y_{i_j}\| \ge |\sum_{i\in n} b_i|$,
and hence, if we pass to a
subsequence of $(y_j)_{j\in\omega}$ having a spreading model, then
this spreading model is equivalent to the unit vector basis of
$\ell_1$.

It remains to show that $\ell_1$ is not isomorphic to an asymptotic
model of $X$.

By the uniform convexity of $\ell_3$ we have:
\begin{gather*}\tag*{($*$)}
\text{for any $\eps>0$ there exists $\lambda<1$ such that}\\ \text{if
$x,y\in B_{\ell_3}$ with $\| x+y\|_{\ell_3}>2\lambda$, then $\|
x-y\|<\eps$\,.}
\end{gather*}

We shall now fix parameters $1>\lambda_1>\lambda_2>\lambda_3>
\lambda_4>\lambda_5>0.9$, $0<\eps_1<\eps_3<\eps_4< \frac14$,
$\delta_4 = 1-\lambda_4$, $\delta_1= 1-\lambda_1$ as follows. We use
$(*)$ to obtain $\lambda_4$ from $\eps_4$, where we require $\eps_4$
(and $\lambda_4$) to satisfy $1-2\delta_4 -2\eps_4>\lambda_5$.
$\lambda_3$ and $\eps_3$ are chosen so that for any normalized basic
sequence $(x_i)_{i\in\omega}$ with a $\lambda_3$-lower $\ell_1$
estimate, if $\| y_i-x_i\| <\eps_3$ for all $i\in\omega$, then $\big(
\frac{y_i}{\|y_i\|}\big)_{i\in\omega}$ admits a $\lambda_4$-lower
$\ell_1$ estimate.
Then choose $\lambda_2$ so that $\lambda_2^3 +\eps_3^3 >1$.
Take $\eps_1>0$ to determine $\lambda_1$ by $(*)$ so
that $1-2\delta_1-\eps_1 >\lambda_2$.
If $\ell_1$ is an asymptotic model of $X$, then, since $X$ is reflexive, by
the proof that $\ell_1$ is
not distortable (cf. \cite{James}), we may assume that $X$ admits a
block basis array $(x_i^n)_{n,i\in\omega}$ which asymptotically
generates $(e_i)_{i\in\omega}$, where $\|\sum_{i\in n}b_i e_i
\|>\lambda_1\sum_{i\in n} |b_i|$ for all scalars $(b_i)_{i\in n}$
not identically zero.
\jump
{\sc Claim:}
For $n\ge 1$ there exists
$K_n\in\omega$ and $i_n\in\omega$ so that if $i\ge i_n$ there exists
$F_i\subs \supp x_i^n$ with $|F_i|\le K_n$ and
$\|x_i^n|_{\omega\setminus F_i}\|<\eps_3$.
\jump\hspace*{\parindent}
To see this, fix $n\ge 1$.
Since $\| e^0+e^n\| > 2\lambda_1$, there exists $k\in\omega$ so
that if $i>k$, then $\| x_k^0+x_i^n\| > 2\lambda_1$.
Let $i>k$ be fixed
and choose disjointly supported $(g_j)_{j\in m} \subs\cG$ so that
\begin{align}\label{eqn1}
\Big( \sum_{j\in m} \big(g_j(x_k^0)
+g_j(x_i^n)\big)^3\Big)^{\frac13} >2\lambda_1\,.
\end{align}
Thus, by our choice of $\eps_1$ using $(*)$,
\begin{align}\label{eqn2}
\Big{\|} \big(g_j(x_k^0)\big)_{j\in m}
-\big(g_j(x_i^n)\big)_{j\in m} \Big{\|}_{\ell_3}< \eps_1\,.
\end{align}
We reorder the $g_j$'s and choose $\bar m \le m$ so that for $j\in
\bar m$, $\supp(g_j)\cap \supp(x_k^0)\neq \emptyset$, and for $j\in
m\setminus \bar m$, $\supp (g_j)\cap \supp (x_k^0)=\emptyset$. {From}
\eqref{eqn1} and the triangle inequality in $\ell_3$,
$\big(\sum_{j\in m} g_j(x_i^n)^3\big)^{1/3} >1-2\delta_1$, and from
\eqref{eqn2} and the choice of $\bar m$ we obtain $( \sum_{j\in
m\setminus \bar m} g_j(x_i^n)^3)^{1/3} <\eps_1$. Thus, by the
triangle inequality,
\begin{align}\label{eqn3}
\Big( \sum_{j\in\bar m} g_j(x_i^n)^3\Big)^{\frac13} >1-2\delta_1
-\eps_1 >\lambda_2\,.
\end{align}
By admissibility restrictions for $j\in\bar m$, $|\supp (g_j)|\le
\max(\supp(x_k^0))$ and thus, since $\bar m \le\max (\supp (x_k^0))$,
$|\bigcup_{j\in \bar m}\supp(g_j)|\le ( \max( \supp(x_k^0)))^2
=:K_n$. Let $F_i=\bigcup_{j\in\bar m} ( \supp(g_j)\cap\supp(x_i^n))$,
so $|F_i|\le K_n$. By \eqref{eqn3}, $1= \| x_i^n\| > \big(\lambda_2^3
+\|x_i^n |_{\omega\setminus F_i}\|^3\big)^{1/3}$ and so, by our
choice of $\lambda_2^3+\eps_3^3 >1$ we obtain
$\|x_i^n|_{\omega\setminus F_i}\|
< \eps_3$, which proves the claim.

Using the claim for $n\ge 1$, let $y_i^n=\frac{x_i^n
|_{\omega\setminus F_i}} {{\|}x_i^n |_{\omega\setminus F_i}{\|}}$ for
$i> i_n$ and $y_i^n=x_i^n$ for $i\le i_n$. By
Proposition\;\ref{prp:K-basic}, we pass to a subarray asymptotically
generating $(f_i)_{i\in\omega}$. By our choice of $\eps_3$ and the
claim, for all not identically zero scalars $(b_i)_{1\le i\le n}$,
$\|\sum_{i=1}^n b_i f_i\| > \lambda_4\cdot\sum_{i=1}^n |b_i|$. Since
$|\supp (y_i^n)|\le K_n$ for $n\ge 1$, by passing to another subarray
we may assume that for $n\ge 1$ there exists $x^n\in c_{00}$ so that
if $i\ge n$ and $y_i^n= (0,\ldots,0,a_1^n,,\ldots,0, a_2^n,0,
\ldots,0,a_{p_n}^n,0,\ldots)$ where the $a_k^n$'s are the non-zero
coordinates of $x_i^n$, then $x^n=(
a_1^n,\ldots,a_{p_n}^n,0,0,\ldots)$. Of course, $p_n\le K_n$. In
short, the $y_i^n$'s are an identically distributed normalized block
basis of $(u_j)_{j\in\omega}$ and $(v_j)_{j\in\omega}$, i.e., in both
$X$ and $Y$ norms. This is done by passing to a subsequence in each
row, iteratively, so that the distributions converge to that of
$x^n$. We then diagonalize. This array still asymptotically generates
$(f_i)_{i\in\omega}$. Of course, we lost our $0^\te$ row, so, let us
relabel every thing as $(x_i^n)_{n,i\in\omega}$ asymptotically
generating $(f_i)_{i\in\omega}$ with the $\lambda_4$-lower $\ell_1$
estimates and the fact the $x_i^n$ equals $x^n$ in distribution for
$i\ge n$. And our old $K_n$ becomes $K_{n-1}$ in the new labeling.
\jump {From} this point on we work in $Y$ (when computing $\|
\sum_{i\in m} b_i x_{i_n}^n\|$ for $i_0$ large, the $X$ and $Y$ norms
coincide). For $x=(a_0,\ldots,a_{n-1},0,0,\ldots)\in c_{00}$, let
$x^*:=(a_{\pi(0)},\ldots,$ $a_{\pi(n-1)},0,0,\ldots)$, where $\pi$ is
a permutation of $n$ such that $|a_{\pi(0)}|\ge\ldots\ge
|a_{\pi(n-1)}|$. By passing to a subsequence of the rows (the new
array still asymptotically generates $\ell_1$ with lower estimate
$\lambda_4$; indeed, it generates a subsequence of
$(f_i)_{i\in\omega}$) we may assume that ${x^n}^* \to x\in c_0$
coordinatewise, where $x=(a_0,a_1, \ldots)$ with $|a_0|\ge
|a_1|\ge\ldots$. Also, since $Y$ is reflexive, $x\in Y$ and $\|
x\|_Y\le 1$. Choose $p\in\omega$ so that
$\|(a_p,a_{p+1},\ldots)\|_Y<\eps_4$; choose $M\in\omega$ so that
$\frac{1}{\sqrt{m_M}}K_0^2 <\eps_4$ (recall that $K_0$ is the
cardinality of the support of $x^0$); and further choose $N>8K_0 M$
so that $(pN)^{1/3}<\frac{N}{8}$.

We next choose $\gamma_n\downarrow 0$ with $\sum_{n\in\omega\setminus
\{0\}} \gamma_n <1$. For each $n\in\omega$ choose
$\bar\gamma_{n+1}>0$ so that if $g= \frac{1_E}{\sqrt{m_i}}$ is a term
of some $f\in \cF$ with the property that $|g(z)| \ge \gamma_n$ for
some $\|z\|_Y \le 1$ with $|\supp (z)| \le K_n$ then $|g(y)|
<\gamma_{n+1}$ whenever $\|y\|_Y \le1$ and $\|y\|_\infty <
\bar\gamma_{n+1}$. By passing to a subsequence of the rows again and
relabeling and not changing the first row of $x_i^0$'s we may assume
that ${x^n}^*|_p=x|_p$ for all positive $n$ (this actually introduces
a slight error which we shall ignore in that it is insignificant to
what follows) and
\begin{align}\label{eqn4}
{x^n}^* = x|_p +x|_{[p,p_n]}
+{x^n}^*|_{(p_n,K_n]}
\end{align}
where $\|{x^n}^*|_{[p_n,K_n]} \|_\infty < \bar \gamma_n$,
the $\|\cdot\|_\infty$ being calculated
relative to the $(u_j)$-coordinates, where
$p<p_1< K_1< p_2< K_2< p_3\ldots$.
Now $\| x_{i_0}^0+\frac{1}{N}\sum_{n=1}^N x_{i_n}^n\|
>2\lambda_4$, provided $i_0<i_1<\ldots <i_N$ are large enough.
We fix these elements and use \eqref{eqn4} to write each $x_{i_n}^n=
x_{i_n}^n(1)+x_{i_n}^n(2)+x_{i_n}^n(3)$, where the three terms are
disjointly supported and each has, respectively, the same
distribution as the three terms in \eqref{eqn4}, thus,
$x_{i_n}^n(2)^* = x|_{[p,p_n]}^*$.
Choose disjointly supported
$(g_k)_{k\in m}\subs\cF$ with
\jump \hspace*{0.25\linewidth}
\begin{equation}\label{eqn5}
\Big( \sum_{k\in m} g_k \big( x_{i_0}^0+\frac{1}{N}\sum_{n=1}^N
x_{i_n}^n\big)^3 \Big)^{\frac13}> 2\lambda_4\,.
\end{equation}
It follows that
\begin{equation*}
\Big( \sum_{k\in m} g_k\big(x_{i_0}^0 \big)^3\Big)^{\frac13}
>1-2\delta_4\,.
\end{equation*}
Write $g_k= \sum_{j\in n_k}
\frac{1_{E_{i_j^k}}}{\sqrt{m\indd}_{\mindd i_j^k}}$ as in the
definition of $\cF$. We shall call
$\frac{1_{E_{i_j^k}}}{\sqrt{m\indd}_{\mindd i_j^k}}$ a {\bf term\/}
of $g_k$. By reordering the $g_k$'s we may assume for some $\bar m\le
m$ that if $k\le\bar m$, then some term $\frac{1_E}{\sqrt{m_j}}$ of
$g_k$ satisfies $|\frac{1_E}{\sqrt{m_j}}(x_{i_0}^0)| \ge
\frac{\eps_4}{K_0}$. In particular, this forces $\bar m\le K_0$ and
$j< M$ and so $n_k <M$ for $k\le \bar m$. If $k\in m\setminus \bar
m$, then for each term $\frac{1_E}{\sqrt{m_j}}$ of $g_k$ we have
$|\frac{1_E}{\sqrt{m_j}}(x_{i_k}^0)| \ge \frac{\eps_4}{K_0}$ and so,
since at most $K_0$ such terms could be non-zero on $x_{i_0}^0$,
\begin{equation}\label{eqn6}
\Big( \sum_{k\in m\setminus \bar m} g_k\big( x_{i_0}^0\big)^3
\Big)^{\frac13}\le \sum_{k\in m\setminus\bar m}\big{|}g_k
(x_{i_0}^0)\big{|} <\frac{\eps_4}{K_0}\cdot K_0=\eps_4\,.
\end{equation}
{From} $(*)$, \eqref{eqn5}      
and our choice of $\lambda_4$, $( \sum_{k\in m}
( g_k(x_{i_0}^0)-g_k(\frac{1}{N}\sum_{n=1}^N x_{i_n}^n))^3
)^{1/3} <\eps_4$, and so, from \eqref{eqn6} and the triangle
inequality in $\ell_3$,
\begin{equation}\label{eqn7}
\bigg( \sum_{k\in m\setminus\bar m} g_k\Big(\frac{1}{N}\sum_{n=1}^N
x_{i_n}^n\Big)^3 \bigg)^{\frac13} <2\eps_4\,.
\end{equation}
Thus, by \eqref{eqn7} and \eqref{eqn5},
\begin{equation}\label{eqn8}
\bigg( \sum_{k\in \bar m} g_k\Big(\frac{1}{N}\sum_{n=1}^N
x_{i_n}^n\Big)^3 \bigg)^{\frac13} >1-2\delta_4 -2\eps_4>\lambda_5\,.
\end{equation}
Now $\bar m\le K_0$ and each $n_k\le M$.
So we have amongst $(g_k)_{k\in\bar m}$ at most $K_0 M$ terms of the form
$\frac{1_E}{\sqrt{m}}$.
We shall show that
\begin{multline}\label{eqn9}
\bigg( \sum_{k\in \bar m} g_k\Big(\frac{1}{N}\sum_{n=1}^N
x_{i_n}^n(1)\Big)^3 \bigg)^{\frac13} + \bigg( \sum\limits_{k\in \bar
m} g_k\Big(\frac{1}{N}\sum_{n=1}^N x_{i_n}^n(2)\Big)^3
\bigg)^{\frac13}\\ + \bigg( \sum_{k\in \bar m} g_k\Big(\frac{1}{N}
\sum_{n=1}^N x_{i_n}^n(3)\Big)^3 \bigg)^{\frac13}< \lambda_5\,,
\end{multline}
which will contradict \eqref{eqn8}. The second term is easiest to
estimate, it is $$\le \frac{1}{N}\sum_{n=1}^{N}\big\| x_{i_n}^n(2)
\big\| <\frac{1}{N}\sum_{n=1}^N \eps_4=\eps_4\,.$$ We next estimate
the third term in \eqref{eqn9}. If for a term
$\frac{1_E}{\sqrt{m_j}}$ of some $g_k$, $k\in\bar m$ we have
$\big{|}\frac{1_E}{\sqrt{m_j}} (x_{i_n}^n(3))\big{|} \ge\gamma_n$,
then $|\frac{1_E}{\sqrt{m_j}} (x_{i_l}^l(3))| \le\gamma_l$ for $l\neq
n$. Thus, $$\Big|\frac{1_E}{\sqrt{m_j}} \Big( \frac{1}{N}
\sum_{n=1}^N x_{i_n}^n(3) \Big)\Big| \le \frac{1}{N}\Big(
1+\sum_{j=1}^N\gamma_j\Big)$$ and therefore the third term in
\eqref{eqn9} is $$\le \frac{1}{N} (K_0 M)\bigg(1+
\sum_{j=1}^{N}\gamma_j \bigg)
< \frac{2 K_0 M}{N}\,.$$
Finally, $\frac{1}{N}\sum_{n=1}^N x_{i_n}^n(1)$ consists of the
vector $\frac{1}{N}x|_p$ repeated $N$ times on disjoint blocks.
Hence, its norm is less than or equal to twice the norm of the vector
in $Y$ which consists of $\frac1N$ repeated $pN$ times. Since
$\sum_{n\in\omega\setminus \{0\}}\frac{1}{\sqrt{m_n}}<1$, this is at
most $\frac{2(pN)^{1/3}}{N}<\frac18$. Thus, the left hand side of
\eqref{eqn9} is $$\le \frac18 +\eps_4 + \frac{2K_0 M}{N} <\frac18
+\frac14 +\frac18 =\frac12 <\lambda_5$$ and we have a contradiction
which completes the proof of Example\;\ref{exa}.
\end{proof}

In summary, asymptotic models generalize spreading models. Certain
positive theorems that one would like to have for spreading models
are just not true. This was one motivation behind the development of
asymptotic structures $\{X\}_n$ in \cite{Maurey.etal}. In that
setting, the theorems are more complete, yet a sacrifice is made in
that certain infinite dimensional structural ties are lost.
Asymptotic models provide a somewhat fuller theory than spreading
models, although some of the same deficiencies remain. They also
provide a context in which some of the long outstanding problems in
spreading models may prove tractable in this new setting (see
Section\;\ref{sec:problems} below for some of these problems). We
believe that the stronger type of convergence one has in strong
asymptotic models, as opposed to the convergence of arrays should
enter into the solution of some of these problems.

\section{Asymptotic Models Under Renormings}\label{sec:renorm}
\setcounter{equation}{0}
\setcounter{thm}{0}

In this section we extend some of the results of
\cite{OS98b}    
to the settings of asymptotic models. Information about the spreading
models of a space $X$ does not usually yield information about the
subspace structure of $X$. For example, every $X\subseteq T$
(Tsirelson's space) has a spreading model 1-equivalent to the unit
vector basis of $\ell_1$, but $T$ does not contain an isomorph of
$\ell_1$ \cite{OS98}. But something can be said if one strengthens
the hypothesis to include all equivalent norms.

\begin{thm}\label{thm:4-1}
{\rm \cite{OS98b}} For every $X$ there exists an equivalent norm
$\Norm\cdot\Norm$ on $X$, so that we have: If $(X,\Norm\cdot\Norm)$
admits a spreading model $(e_n)_{n\in\omega}$ satisfying
\begin{itemize}
\item [a)] $(e_n)_{n\in\omega}$ is 1-equivalent to the unit vector
basis of $c_0$ (or even just $\Norm e_0 + e_1 \Norm =1$, where
$(e_n)_{n\in\omega}$ is generated by a weakly null sequence), then
$X$ contains an isomorph of $c_0$;
\item [b)] $(e_n)_{n\in\omega}$ is 1-equivalent to the unit vector
basis of $\ell_1$ (or even just $\Norm e_0 \pm e_1\Norm =2$), then
$X$ contains an isomorph of $\ell_1$;
\item [c)] $(e_n)_{n\in\omega}$ is such that $\|\sum_{i\in\omega}
a_i e_i\| = \sum_{i\in\omega} a_i$ for all $(a_i)\in c_{00}$ with
$a_i \ge 0$ for $i\in\omega$ (or even just $\Norm e_0 +e_1\Norm =2$),
then $X$ is not reflexive.
\end{itemize}
\end{thm}

We shall develop an asymptotic model version of each part. Part of
our construction will mirror that in \cite{OS98b}, but we need some
new tricks as well. We begin by recalling the construction of the
equivalent norm $\Norm\cdot\Norm$ from \cite{OS98b}.

For $c\in X$ and $x\in X$ define $\|x\|_c : = \big{\|}c\|x\| + x
\big{\|} +  \big{\|}c\|x\| - x \big{\|}$, where $\|\cdot\|$ is the
original norm on $X$. Then $\|x\|_c$ is an equivalent norm on $X$ and
in fact, for all $x\in X$, $2\|x\| \le \|x\|_c \le 2 (1+\|c\|)\|x\|$.
Let $C$ be a countable dense set in $X$ and for $c\in C$ choose
$p_c>0$ so that $\sum_{c\in C} p_c (1+\|c\|) <\infty$. Define for
$x\in X$,
\begin{equation}\label{eqn:4-1}
\Norm x\Norm : = \sum_{c\in C} p_c \|x\|_c\,.
\end{equation}
This is an equivalent norm on $X$. We call $\Norm\cdot\Norm$ the {\bf
asymptotic norm\/} generated by $\|\cdot\|$. We may assume $\Norm
x\Norm \ge \|x\|$.

\begin{thm}\label{thm:4-2}
$X$ contains an isomorph of $c_0$ if there exists a weakly null basic
array $(x_i^n)_{n,i\in\omega} \subseteq X$ generating---in $(X,\Norm
\cdot\Norm)$---an asymptotic model $(e_i)_{i\in\omega}$ which is
1-equivalent to the unit vector basis of $c_0$.
\end{thm}

\begin{lem}\label{lem:4-3}
Let $(x_m)_{m\in\omega}$ and $(y_n)_{n\in\omega}$ be $\Norm \cdot
\Norm$ normalized weakly null sequences in $X$ with
$\lim_{m\to\infty} \lim_{n\to\infty} \Norm x_m + y_n \Norm =1$. Then
there exist integers $k(0) < k(1) <\cdots$ so that setting $a:=
\lim_{m\to\infty} \|x_{k(m)}\|$ and $x'_m = \frac{x_{k(m)}}{\|
x_{k(m)}\|}$, for all $y\in X$ we have
\begin{equation}\label{eqn:4-2}
\lim_{m\to\infty} \lim_{n\to\infty} \|y +x'_m + a^{-1} y_{k(n)}\|
 = \lim_{m\to\infty} \|y + x'_m\|\,.
\end{equation}
\end{lem}

\begin{proof}
By Ramsey's Theorem there exist $k(0) < k(1) <\cdots $ so that for
all $y\in X$ and $\alpha,\beta \in{\mathbb R}$,
$$\lim_{m\to\infty}
\lim_{n\to\infty} \|y + \alpha x_{k(m)} +\beta y_{k(n)} \| \ \mbox{
exists.}$$
To simplify notation we write $(x_m)_{m\in\omega}$ and
$(y_n)_{n\in\omega}$ for $(x_{k(m)})_{m\in\omega}$ and
$(y_{k(n)})_{n\in\omega}$ and thus $a: = \lim_{m\to\infty} \|x_m\|$.
Now
\begin{equation*}
\begin{split}
1 & = \lim_{m\to\infty} \lim_{n\to\infty} \Norm x_m +y_n\Norm
= \lim_{m\to\infty} \lim_{n\to\infty} \sum_{c\in C} p_c \|x_m + y_n\|_c\\
& = \lim_{m\to\infty} \sum_{c\in C} p_c \|x_m\|_c\ .
\end{split}
\end{equation*}
Thus
\begin{equation}\label{eqn:4-3}
1 = \sum_{c\in C} p_c \Big( \lim_{m\to\infty}\lim_{n\to\infty} \|x_m
+ y_n\|_c\Big) = \sum_{c\in C} p_c \Big( \lim_{m\to\infty}
\|x_m\|_c\Big)\,.
\end{equation}
Since $y_n \to0$ weakly, $\lim_{m\to\infty} \lim_{n\to\infty}
\|x_m+y_n\|_c \ge \lim_{m\to\infty} \|x_m\|_c$ for all $c\in C$.
{From} this and \eqref{eqn:4-3} we get
\begin{equation}\label{eqn:4-4}
\lim_{m\to\infty} \lim_{n\to\infty} \|x_m + y_n \|_c
= \lim_{m\to\infty} \|x_m\|_c \text{ for all } c\in X
\text{ (since $C$ is dense).}
\end{equation}
Letting $c=0$, this yields
\begin{equation*}
\lim_{m\to\infty} \lim_{n\to\infty} \|x_m + y_n\| =a\ .
\end{equation*}
Thus, for all $y\in X$,
\begin{equation}\label{eqn:4-5}
\begin{split}
&\lim_{m\to\infty} \lim_{n\to\infty} \big[ \|ay + x_m + y_n\| + \|
-ay +x_m + y_n\|\big] \\ &\qquad\qquad = \lim_{m\to\infty} \big[ \|ay
+ x_m\| + \|-ay + x_m\|\big]\ .
\end{split}
\end{equation}
Again, since $(y_n)_{n\in\omega}$ is weakly null,
\begin{equation*}
\begin{split}
\lim_{n\to\infty} \|ay + x_m + y_n\|
& \ge \|ay + x_m\| \ \mbox{and}\\
\lim_{n\to\infty} \|-ay + x_m + y_n\|
& \ge \|-ay+ x_m\| \  \mbox{ for }\ m\in\omega\ .
\end{split}
\end{equation*}
Thus, by \eqref{eqn:4-5}, for all $y\in X$ we have
$$\lim_{m\to\infty} \lim_{n\to\infty} \|ay + x_m + y_n\| =
\lim_{m\to\infty} \|ay + x_m\|$$ which completes the proof.
\end{proof}

Note that it follows from Lemma\;\ref{lem:4-3} that if
$\lim_{m\to\infty} \lim_{n\to\infty}  \Norm x_m \pm y_n\Norm =1$,
then for all $y\in X$ we can obtain
\begin{equation}\label{eqn:4-6}
\lim_{m\to\infty} \lim_{n\to\infty} \|y \pm x'_m \pm a^{-1} y_{k(n)}\|
= \lim_{m\to\infty} \|y\pm x'_m\|
\end{equation}
for all choices of sign (keeping the sign of $x'_m$ the same on both
sides of \eqref{eqn:4-6}).

\begin{proof}[Proof of Theorem\;\ref{thm:4-2}]
By passing to a subarray of $(x_i^n)_{n,i\in\omega}$ we may assume
that for each $n\in\omega$ we have $\lim_{i\to\infty} \|x_i^n\| =
a_n$ (for some $a_n$). Let $\eps_n\downarrow 0$ with
$\sum_{n\in\omega} \eps_n <\infty$. By passing to a subsequence of
the rows we may assume that for all $n$, $a_n\to a>0$,
$|\frac1{a_n}-\frac1{a} |<\frac{\eps_n}{3}$ and $a_n > \frac{a}{2}$.
In addition we may assume that for all $y \in X$,
$\alpha,\beta\in{\mathbb R}$ and $i,j\in\omega$, $$\lim_{m\to\infty}
\lim_{n\to\infty} \| y + \alpha x_m^i + \beta x_n^j\| $$ exists, and
moreover, by Lemma\;\ref{lem:4-3} (actually \eqref{eqn:4-6}) we may
assume that for $y\in X$ and $p,q\in\omega$ with $p<q$ we have
$$\lim_{i\to\infty} \lim_{j\to\infty} \Big\| y\pm\frac{x_i^p}{a_p}
\pm \frac{x_j^q}{a_p}\Big\| = \lim_{i\to\infty} \Big\| y \pm
\frac{x_i^p}{a_p}\Big\|\,.$$ Hence, from the triangle inequality
using $|\frac1a - \frac1{a_n}| < \frac{\eps_n}{3}$ we get
\begin{equation}\label{eqn:4-7}
\lim_{i\to\infty} \lim_{j\to\infty} \Big\| y \pm \frac{x_i^p}a \pm
\frac{x_j^q}a\Big\|
< \lim_{i\to\infty} \Big\| y + \frac{x_i^p}a\Big\| +\eps_i\ .
\end{equation}

By passing to another subarray and setting $x_i = \frac{x_i^i}{a}$
for $i\in\omega$ we may assume that for all $m\in\omega$ and $y\in 2m
B_{\langle x_i\rangle_{i\in m}}$,
\begin{equation}\label{eqn:4-8}
\|y \pm x_m \pm x_{m+1}\| < \| y \pm x_m\| + 2\eps_m\,.
\end{equation}

This is accomplished using \eqref{eqn:4-7}. If $i$ is large enough
and $i<j$, then $\big{\|} \frac{x_i^0}a \pm \frac{x_j^1}a\big{\|}
< \big{\|}\frac{x_0^i}a\big{\|} +\eps_0$.
This fixes $i$ and $x_0 =  \frac{x_i^0}a$ (under relabeling) and then
we increase $j$ large enough so that  for $j<k$ and
$y\in 2B_{\langle x_0\rangle}$,
$$\Big\| y \pm \frac{x_j^1}a \pm \frac{x_k^2}a\Big\|
< \Big\| y + \frac{x_j^1}a\Big\| + \eps_1\ .$$
This fixes $j$ and $x_1 = \frac{x_j^1}a$ (under relabeling) and so on.

We claim that $$\sup \biggl\{ \Big\| \sum_{i\in m} \pm x_i\Big\| :
\mbox{ all choices of } \pm \biggr\} < \infty\,,$$ which will yield
the theorem ($(x_i)_{i\in\omega}$ is then equivalent to the unit
vector basis of $c_0$). Indeed, from \eqref{eqn:4-7} we get
\begin{equation*}
\begin{split}
\Big\| \sum_{i\in m} \pm x_i\Big\| & \le \Big\| \sum_{i\in m-1} \pm
x_i \Big\| + 2\eps_{m-2}\\ & \le \cdots \le \|x_0\| +
\sum_{m\in\omega} 2\eps_m < \infty\,.
\end{split}
\end{equation*}\vspace*{-4ex}

\end{proof}
\smallskip

\begin{rmk}\label{rmk:4-4}
In the proof of Theorem\;\ref{thm:4-2} we only used
$\lim_{m\to\infty} \lim_{n\to\infty} \Norm x_m^p \pm x_n^q\Norm =1$
for all $p<q$. In other words $\Norm e_p \pm e_q\Norm =1$ for $p\ne
q$. In the case of spreading models (Theorem\;\ref{thm:4-1}(a)) one
only needs $\Norm e_p + e_q\Norm =1$ for $p\ne q$. We do not know if
this is sufficient to obtain $c_0$ inside $X$ for asymptotic models.
\end{rmk}

The proof of Theorem\;\ref{thm:4-2} was the most similar to the
spreading model analogue of the three results we present in this
section. Our next proof is more difficult.

\begin{thm}\label{thm:4-5}
For every separable infinite dimensional Banach space $X$, there
exists an equivalent norm $\NORM \cdot\NORM$ on $X$ with the
following property. If there exist $\NORM\cdot\NORM$-normalized basic
sequences $(x_m)_{m\in\omega}$ and $(y_n)_{n\in\omega}$ with
$\lim_{m\to\infty}\lim_{n\to\infty} \NORM x_m + y_n\NORM =2$, then
$X$ is not reflexive.
\end{thm}

\begin{cor}\label{cor:4-6}
$X$ is reflexive if and only if there exists an equivalent norm
$\NORM \cdot\NORM$ on $X$ such that if $(e_n)_{n\in\omega}$ is an
asymptotic model of $(X,\NORM \cdot\NORM )$, then $\NORM e_0 +
e_1\NORM <2$.
\end{cor}

\begin{proof}[Proof of Theorem\;\ref{thm:4-5}]
We first construct the norm $\NORM \cdot\NORM $ on $X$. We begin by
assuming that $X= \langle x_0\rangle \oplus_\infty Y$ where $Y$ is a
subspace of a Banach space with a bimonotone normalized basis $(d_i)$
and we let $)\!)\cdot(\!($ be the inherited norm on $Y$. We assume
the norm $\|\cdot\|$ on $X$ is given as follows. If $x = ax_0 +y\in
X$ with $a\in {\mathbb R}$ and $y\in Y$, then $\|x\| = \max
(|a|,)\!)y(\!( + \sum_{i\in\omega} |y(i)|2^{-i})$ if $y =
\sum_{i\in\omega} y(i) d_i$. We have the following:
\begin{eqnarray}
&&\mbox{Let $(x_m)_{m\in\omega}$ and $(y_n)_{n\in\omega}$ be weakly
null $\|\cdot\|$ normalized sequences in $X$.} \label{eqn:4-9} \\
&&\mbox{Let $\alpha +\beta =1$, $\alpha,\beta>0$ and $\alpha \ne
\frac12$.}\nonumber\\ &&\mbox{Then
$\lim_{m\to\infty}\lim_{n\to\infty} \|x_0 +\frac12 x_m + \frac12
y_n\|=1$ while}\nonumber\\ &&\displaystyle \lim_{m\to\infty} \big\|
\alpha x_0 + \tfrac12 x_m\big\| + \lim_{n\to\infty} \big\| \beta x_0
+ \tfrac12 y_n\big\| \nonumber\\ &&\qquad = \max \big(
\alpha,\tfrac12\big) + \max \big(\beta,\tfrac12\big) = \tfrac12 +
\max (\alpha,\beta) >1\,.\nonumber
\end{eqnarray}
\begin{eqnarray}
&&\mbox{Let $y\in Y$, $y\ne 0$ and let $(x_m)_{m\in\omega}$ be a
$\|\cdot\|$-normalized weakly null }\label{eqn:4-10}\\
&&\mbox{sequence in $X$. Then, presuming the limit exists,}\nonumber\\
&&\qquad \displaystyle
\lim_{n\to\infty} \|y +x_n\| \ge 1 + \sum_{i\in\omega} 2^{-i} |y(i)|
>1\,.\nonumber
\end{eqnarray}

Let $\Norm\cdot\Norm$ be the asymptotic norm on $X$ generated by
$\|\cdot\|$ (see \eqref{eqn:4-1} above), and let $\NORM\cdot\NORM$ be
the equivalent asymptotic norm on $X$ generated by
$\Norm\cdot\Norm$\,.
\renewcommand \qedsymbol {}
\end{proof}
\renewcommand \qedsymbol {$\boldsymbol{\dashv}$}

Before proceeding we present a lemma. The lemma is valid in any
$(X,\|\cdot\|)$, not just in our space above.

\begin{lem}\label{lem:4-7}
Let $\Norm\cdot\Norm$ be the equivalent asymptotic norm on
$(X,\|\cdot\|)$ generated by $\|\cdot\|$ as in \eqref{eqn:4-1}. Let
$(x_m)_{m\in\omega}$ and $(y_n)_{n\in\omega}$ be
$\Norm\cdot\Norm$-normalized sequences in $X$.
\begin{itemize}
\item[a)] If\/ $\lim_{m\to\infty}\lim_{n\to\infty} \Norm x_m +
y_n\Norm=2$, then there exist integers $k(0) < k(1) <\cdots$ so that
setting $x'_m = \frac{x_{k(m)}}{\|x_{k(m)}\|}$ and $y'_n =
\frac{y_{k(n)}}{\|y_{k(n)}\|}$, then for all $y\in Y$ and
$\beta_1,\beta_2 \ge 0$ (not both~$0$) we have
\begin{align}
&\lim_{m\to\infty} \lim_{n\to\infty} \|y +\beta_1  x'_m + \beta_2
y'_n\| \tag{11a}\label{eqn:4-11a}\\ &\qquad = \lim_{m\to\infty}
\Big\| \frac{\beta_1}{\beta_1+\beta_2} y + \beta_1 x'_m\Big\| +
\lim_{n\to\infty}  \Big\| \frac{\beta_2}{\beta_1+\beta_2} y + \beta_2
y'_n\Big\|\,. \nonumber
\end{align}
\item[b)]
If\/ $\lim_{m\to\infty} \lim_{n\to\infty} \Norm x_m\pm y_n\Norm=2$,
then there exist integers $k(0) < k(1) < \cdots$ so that setting
$x'_m = \frac{x_{k(m)}}{\|x_{k(m)}\|}$ and $y'_n =
\frac{y_{k(n)}}{\|y_{k(n)}\|}$, then for all $y\in X$,
$\beta_1,\beta_2 \in{\mathbb R}$ (not both~$0$) we have
\begin{align}
&\lim_{m\to\infty} \lim_{n\to\infty} \|y +\beta_1  x'_m + \beta_2
y'_n\| \tag{11b}\label{eqn:4-11b}\\ &\qquad = \lim_{m\to\infty}
\Big\| \frac{|\beta_1|}{|\beta_1|+|\beta_2|} y + \beta_1 x'_m\Big\| +
\lim_{n\to\infty}  \Big\| \frac{|\beta_2|}{|\beta_1|+|\beta_2|} y +
\beta_2 y'_n\Big\|\,. \nonumber
\end{align}
\addtocounter{equation}{1}
\end{itemize}
\end{lem}

\begin{proof}
Again by Ramsey's Theorem we can find $k(0) < k(1) <\cdots$ so that
relabeling $x_{k(m)} = x_m$ and $y_{k(n)} = y_n$,
$\lim_{m\to\infty}\lim_{n\to\infty} \|y +\alpha x_m +\beta y_n\|$
exists for all $y\in X$ and $\alpha,\beta \in {\mathbb R}$. Let $a =
\lim_{m\to\infty} \|x_m\|$, $b= \lim_{n\to\infty} \|y_n\|$ and let
$x'_m = \frac{x_m}{a}$, $y'_n = \frac{y_n}{b}$. We will prove the
conclusion of the lemma for these sequences which will yield the
lemma.

{\bf a)} We first suppose that $\beta_1 + \beta_2=1$. Set
$\bar\beta_1 = \frac{\beta_1}{a}$, $\bar\beta_2 = \frac{\beta_2}{b}$.
{From} our hypothesis, $$\lim_{m\to\infty} \lim_{n\to\infty} \Norm
\bar \beta_1 x_m + \bar\beta_2 y_n\Norm = \bar\beta_1 + \bar
\beta_2\,.$$

{From} the definition of $\Norm\cdot\Norm$ and the triangle
inequality in each $\|\cdot\|_c$ we obtain
\begin{equation}\label{eqn:4-12}
\mbox{for $c\in C$, }\ \lim_{m\to\infty} \lim_{n\to\infty} \|\bar
\beta_1 x_m + \bar \beta_2 y_n\|_c = \lim_{m\to\infty} \|\bar\beta_1
x_m\|_c + \lim_{n\to\infty} \|\bar\beta_2 y_n\|_c\,.
\end{equation}
By the density of $C$ in $X$ this holds for all $c\in X$.

Setting $c=0$ in \eqref{eqn:4-12} yields
\begin{equation}\label{eqn:4-13}
\lim_{m\to\infty} \lim_{n\to\infty} \|\beta_1 x'_m + \beta_2 y'_n\| =
\beta_1 +\beta_2 = 1\,.
\end{equation}

{From} \eqref{eqn:4-12}, using \eqref{eqn:4-13}, for all $c\in X$,
\begin{align}
&\lim_{m\to\infty} \lim_{n\to\infty} \Big[ \big{\|}c+\beta_1 x'_m +
\beta_2 y'_n\big{\|} + \big{\|}c - (\beta_1 x'_m + \beta_2
y'_n)\big{\|}\Big] \label{eqn:4-14}\\ &\qquad\qquad =
\lim_{m\to\infty} \Big[ \big{\|}\beta_1 c +\beta_1 x'_m\big{\|} +
\big{\|}\beta_1 c - \beta_1 x'_m\big{\|}\Big] \nonumber\\ &\qquad
\qquad\qquad + \lim_{n\to\infty} \Big[ \big{\|}\beta_2 c+\beta_2
y'_n\big{\|} + \big{\|} \beta_2 c - \beta_2 y'_n\big{\|}
\Big]\,.\nonumber
\end{align}

{From} \eqref{eqn:4-14} and the triangle inequality we obtain
\eqref{eqn:4-11a} in the case $\beta_1 +\beta_2 =1$. To get the
general case from this we note that for $y\in X$,
$\beta_1,\beta_2\in{\mathbb R}$ (not both~$0$) we have
\begin{align*}
&\lim_{m\to\infty} \lim_{n\to\infty} \Big\| \frac{y}{\beta_1+\beta_2}
+ \frac{\beta_1}{\beta_1+\beta_2} x'_m +
\frac{\beta_2}{\beta_1+\beta_2} y'_n\Big\|\\ &\qquad\qquad =
\lim_{m\to\infty} \Big\| \frac{\beta_1}{\beta_1 +\beta_2}
\Big(\frac{y}{\beta_1+\beta_2}\Big) + \frac{\beta_1}{\beta_1+\beta_2}
x'_m \Big\|\\ &\qquad\qquad\qquad + \lim_{n\to\infty} \Big\|
\frac{\beta_2}{\beta_1+\beta_2} \Big(\frac{y}{\beta_1+\beta_2}\Big) +
\frac{\beta_2}{\beta_1+\beta_2} y'_n\Big\|\,
\end{align*}
and \eqref{eqn:4-11a} follows by multiplying by $\beta_1 +\beta_2$.

{\bf b)} We continue the argument from a). As in that case we may
assume that $|\beta_1| + |\beta_2| =\nolinebreak 1$. The case
$\beta_1 ,\beta_2 \le 0$ is covered by a) using $$\lim_{m\to\infty}
\lim_{n\to\infty} \| y +\beta_1 x'_m +\beta_2 y'_n\| =
\lim_{m\to\infty} \lim_{n\to\infty} \| -y-\beta_1 x'_m-\beta_2
y'_n\|\ .$$ Similarly, the only case left to consider is $\beta_1>0$
and $\beta_2<0$. We prefer to take $\beta_1,\beta_2 >0$, $\beta_1
+\beta_2=1$ and work with ``$\beta_1x'_m - \beta_2 y'_n$''. As in a),
we obtain from the hypothesis for $c\in X$, $$\lim_{m\to\infty}
\lim_{n\to\infty} \|\beta_1 x'_m-\beta_2 y'_n\|_c = \lim_{m\to\infty}
\|\beta_1 x'_m\|_c + \lim_{n\to\infty} \|\beta_2 y'_n\|\,.$$ Thus,
for $y\in X$ we get
\begin{align*}
&\lim_{m\to\infty} \lim_{n\to\infty} \Big[ \big{\|}y + (\beta_1 x'_m
-\beta_2 y'_n)\big{\|} + \big{\|} y- (\beta_1 x'_m - \beta_2
y'_n)\big{\|} \Big]\\ &\qquad\qquad = \lim_{m\to\infty} \Big[
\big{\|}\beta_1 y+\beta_1 x'_m\big{\|} + \big{\|} \beta_1 y - \beta_1
x'_m\big{\|} \Big]\\ &\qquad\qquad\qquad + \lim_{n\to\infty} \Big[
\big{\|} \beta_2 y - \beta_2 y'_n\big{\|} + \big{\|}\beta_2 y +
\beta_2 y'_n\big{\|}\Big]\,.
\end{align*}
Again from the triangle inequality we obtain \eqref{eqn:4-11b} in
this case.
\end{proof}

{\it We return to the proof of Theorem\;\ref{thm:4-5}}. Suppose that
$(x_m)_{m\in\omega}$ and $(y_n)_{n\in\omega}$ are $\NORM\cdot\NORM$
normalized basic sequences in $X$ with $\lim_{m\to\infty}
\lim_{n\to\infty} \NORM x_m+y_n \NORM =2$. Assume towards a
contradiction that $X$ is reflexive. Then $(x_m)_{m\in\omega}$ and
$(y_n)_{n\in\omega}$ are both weakly null. We may assume that
$\lim_{m\to\infty}\lim_{n\to\infty} \|y +\alpha x_m +\beta y_n\|$
exists for all $y\in X$, $\alpha,\beta \in{\mathbb R}$ (and  for all
of the norms we have constructed). By Lemma\;\ref{lem:4-7} we may
also assume that setting $x'_m = \frac{x_m}{\Norm x_m\Norm}$ and
$y'_n = \frac{y_n}{\Norm y_n\Norm}$, for $y\in X$ and $\alpha ,\beta
\ge0$ (not both~$0$) we have
\begin{align*}
&\lim_{m\to\infty} \lim_{n\to\infty} \Norm y+\alpha x'_m +\beta
y'_n\Norm \\ &\qquad = \lim_{m\to\infty} \BigNorm
\frac{\alpha}{\alpha+\beta} y +\alpha x'_m\BigNorm +\lim_{n\to\infty}
\BigNorm \frac{\beta}{\alpha+\beta} y + \beta y'_n\BigNorm\,.
\end{align*}
Thus,
\begin{align*}
&\lim_{m\to\infty} \lim_{n\to\infty}  \sum_{c\in C} p_c \bigg[
\Big{\|}c\| y + \alpha x'_m +\beta y'_n\| + y + \alpha x'_m +\beta
y'_n\Big{\|}\\ &\qquad\qquad + \Big{\|}c\| y + \alpha x'_n +\beta
y'_n\| - (y +\alpha x'_m +\beta y'_n)\Big{\|}\bigg]\\ &\qquad =
\lim_{m\to \infty} \sum_{c\in C} p_c \bigg[ \Big{\|}c\big{\|}
\frac{\alpha}{\alpha+\beta} y +\alpha x'_m \big{\|} +
\frac{\alpha}{\alpha+\beta} y+\alpha x'_m\Big{\|} \\
&\qquad\qquad\qquad + \Big{\|}c\big{\|} \frac{\alpha}{\alpha+\beta} y
+ \alpha x'_m \big{\|}-\Big(\frac{\alpha}{\alpha +\beta} y +\alpha
x'_m\Big)\Big{\|}\bigg]\\ &\qquad\qquad +\lim_{n\to\infty} \sum_{c\in
C} p_c \bigg[ \Big{\|}c\big{\|} \frac{\beta}{\alpha+\beta} y +\beta
y'_n \big{\|} + \frac{\beta}{\alpha +\beta} y+\beta y'_n\Big{\|}\\
&\qquad\qquad\qquad + \Big{\|}c\big{\|} \frac{\beta}{\alpha+\beta} y
+ \beta y'_n \big{\|} - \Big( \frac{\beta}{\alpha+\beta} y +\beta
y'_n\Big)\Big{\|}\bigg]\ .
\end{align*}
{From} this and the triangle inequality we have for all $c\in X$,
$y\in X$ and $\alpha,\beta\ge 0$ (not both~$0$) that
\begin{align}
&\lim_{m\to\infty} \lim_{n\to\infty}\big{\|}c\| y +\alpha x'_m +\beta
y'_n \| + y + \alpha x'_m + \beta y'_n\big{\|} \label{eqn:4-15}\\
&\qquad = \lim_{m\to\infty} \Big{\|}c\big{\|}
\frac{\alpha}{\alpha+\beta} y + \alpha x'_m \big{\|} +
\frac{\alpha}{\alpha+\beta} y +\alpha x'_m\Big{\|}\nonumber\\
&\qquad\qquad + \lim_{n\to\infty} \Big{\|}c\big{\|}
\frac{\beta}{\alpha+\beta} y + \beta y'_n \big{\|} +
\frac{\beta}{\alpha+\beta} y + \beta y'_n\Big{\|}\,.\nonumber
\end{align}
Setting $c=y =0$ in \eqref{eqn:4-15} yields
\begin{equation}\label{eqn:4-16}
\lim_{m\to\infty} \lim_{n\to\infty} \|\alpha x'_m +\beta y'_n\| =
\alpha a +\beta b\,,
\end{equation}
where $a= \lim_m \|x'_m\|$ and $b= \lim_n \|y'_n\|$. Let $x''_m =
\frac{x'_m}{\|x'_m\|}$ and $y''_n = \frac{y'_n}{\|y'_n\|}$. Then
\begin{equation}\label{eqn:4-17}
\lim_{m\to\infty} \lim_{n\to\infty} \|\alpha x''_m + \beta y''_n\| =
\alpha +\beta\,.
\end{equation}

Letting $y=0$ and replacing $c$ by $\frac{c}{\alpha+\beta}$ in
\eqref{eqn:4-15}, using \eqref{eqn:4-17}, we have
\begin{align}
&\lim_{m\to\infty} \lim_{n\to\infty} \|c+\alpha x''_m + \beta y''_n\|
\label{eqn:4-18} \\ &\qquad = \lim_{m\to\infty} \Big\|
\frac{\alpha}{\alpha +\beta} c+\alpha x''_m\Big\| + \lim_{n\to\infty}
\Big\| \frac{\beta}{\alpha+\beta} c+\beta y''_n\Big\|\,. \nonumber
\end{align}

We claim that $a=b$. Indeed, let us assume $a\neq b$. By
\eqref{eqn:4-17} we get $\lim_{m\to\infty} \lim_{n\to\infty}
\|\frac12 x''_m + \frac12 y''_n\| =1$ and further we have
$\lim_{m\to\infty} \lim_{n\to\infty} \|x_0 +\frac12 x''_m +\frac12
y''_n\| =1$, see \eqref{eqn:4-9}. But from \eqref{eqn:4-15}, taking
$y= x_0$ and $c= 0$, we get
\begin{align*}
&\lim_{m\to\infty} \lim_{n\to\infty} \Big\| x_0 +\tfrac12 x''_m
+\tfrac12 y''_n \Big\| = \lim_{m\to\infty} \lim_{n\to\infty} \Big\|
x_0 +\tfrac1{2a} x'_m +\tfrac1{2b} y'_n\Big\| \\ &\qquad =
\lim_{m\to\infty} \left\| \frac{\frac1{2a}}{\frac1{2a} +\frac1{2b}}
x_0 + \tfrac12 x''_m \right\| + \lim_{n\to\infty} \left\|
\frac{\frac1{2b}}{\frac1{2a} +\frac1{2b}} x_0 + \tfrac12 y''_n
\right\| \\ &\qquad\qquad
> 1\ \mbox{ (for $a\ne b$) using \eqref{eqn:4-9}.}
\end{align*}

{From} \eqref{eqn:4-15} we obtain for all $c,y\in X$ and
$\alpha,\beta>0$ (not both~$0$), by replacing $\alpha,\beta$ by
$\frac{\alpha}a$ and $\frac{\beta}a$ since $\frac{\beta}b = \frac{\beta}a$,
\begin{align}
&\lim_{m\to\infty} \lim_{n\to\infty} \Big{\|}c\| y + \alpha x''_m +
\beta y''_n \| + y+\alpha x''_m +\beta y''_n\Big{\|}
\label{eqn:4-19}\\ &\qquad = \lim_{m\to\infty} \Big{\|}c\big{\|}
\frac{\alpha}{\alpha+\beta} y + \alpha x''_m \big{\|} +
\frac{\alpha}{\alpha +\beta} y + \alpha x''_m \Big{\|}\nonumber\\
&\qquad\qquad = \lim_{n\to\infty} \Big{\|}c\big{\|}
\frac{\beta}{\alpha+\beta} y + \beta y''_n \big{\|} +
\frac{\beta}{\alpha +\beta} y + \beta y''_n \Big{\|}\,. \nonumber
\end{align}

Next, we wish to show that $(x''_m)_{m\in\omega}$ and
$(y''_n)_{n\in\omega}$ generate the same type over $Y$, i.e., if
$y\in Y$, $\delta : = \lim_{m\to\infty} \|y + x''_m \|$ and $\gamma :
= \lim_{n\to\infty} \|y + y''_m \|$, then $\delta =\gamma$. Clearly,
$\delta =\gamma =1$ if $y=0$, so assume $y\ne 0$ and
$\delta\ne\gamma$. Let $\alpha +\beta =1$.
Now from \eqref{eqn:4-18} we get
\begin{align*}
\lim_{m\to\infty} \lim_{n\to\infty} \|y+ \alpha x''_m +\beta y''_n\|
& = \lim_{m\to\infty} \|\alpha y +\alpha x''_m\| + \lim_{n\to\infty}
\|\beta y + \beta y''_n\| \\ & = \hspace*{1ex}\alpha \delta + \beta
\gamma\,.
\end{align*}

Thus, from \eqref{eqn:4-19} we get for $c\in X$, $\alpha +\beta=1$
and $\alpha,\beta \ge 0$,
\begin{align}
&\lim_{m\to\infty} \lim_{n\to\infty} \|c (\alpha \delta +\beta\gamma)
+ y + \alpha x''_m +\beta y''_n\| \label{eqn:4-20}\\ &\qquad =
\lim_{m\to\infty} \| (\alpha\delta) c + \alpha y +\alpha x''_m \| +
\lim_{n\to\infty} \| (\beta\gamma) c +\beta y + \beta y''_n\|\,.
\nonumber
\end{align}

Let $\alpha =\beta = \frac12$ and $c =
\frac{-1}{\frac{\gamma}2+\frac{\delta}2} y = \frac{-2y}{\delta
+\gamma}$. Using this in \eqref{eqn:4-20}, from \eqref{eqn:4-17} we
have
\begin{align}
\lim_{m\to\infty} \lim_{n\to\infty} \Big\|\frac12 x''_m +\frac12
y''_n\Big\| = 1 = & \lim_{m\to\infty} \Big\| \Big(\frac12 -
\frac{\delta}{\delta+\gamma}\Big) y + \frac12 x''_m
\Big\|\label{eqn:4-21} \\ &\qquad + \lim_{n\to\infty} \Big\| \Big(
\frac12 - \frac{\gamma}{\delta+\gamma}\Big) y + \frac12 y''_n\Big\|
\nonumber
\end{align}
and since $\delta\ne \gamma$, both coefficients of $y\in Y$ on the
right side of \eqref{eqn:4-21} are nonzero. Therefore, by
\eqref{eqn:4-10}, the right side exceeds~1, a contradiction.

It follows that $\lim_{m\to\infty} \lim_{n\to\infty} \|x''_m +
x''_n\|=2$ and moreover, $(x''_n)_{n\in\omega}$ can be substituted
for $(y''_n)_{n\in\omega}$ in our above equations. So, we are in the
same situation as the proof of Theorem\;4.1\,c) in \cite{OS98b} and
it follows that for some subsequence $(x''_{n_i})_{i\in\omega}$,
$$\bigg\| \sum_{i\in \omega} a_i x''_{n_i}\bigg\| > \frac12 \ \mbox{
if }\ (a_i)_{i\in\omega} \subseteq [0,\infty)\ ,\qquad
\sum_{i\in\omega} a_i =1\,.$$ Hence, $(x''_{n_i})_{i\in\omega}$ is
not weakly null and $X$ is not reflexive, which completes the proof
of Theorem\;\ref{thm:4-5}.\hfill\qedsymbol

\begin{thm}\label{thm:4-8}
Let $X$ have a basis $(b_i)_{i\in\omega}$. There exists an equivalent
norm $\NORM \cdot\NORM$ on $X$ so that if $(X,\NORM\cdot\NORM)$
admits $\NORM\cdot\NORM$ normalized block bases of
$(b_i)_{i\in\omega}$, say $(x_m)_{m\in\omega}$ and
$(y_n)_{n\in\omega}$, satisfying $\lim_{m\to\infty} \lim_{n\to\infty}
\NORM x_m\pm y_n\NORM =2$, then $X$ contains an isomorph of $\ell_1$.
\end{thm}

\begin{cor}\label{cor:4-9}
If $X$ has a basis and does not contain an isomorph of $\ell_1$, then
$X$ can be given an equivalent norm so that if $(e_n)_{n\in\omega}$
is any asymptotic model generated by a block basic array, then
$(e_n)_{n\in\omega}$ is not $1$-equivalent to the unit vector basis
of $\ell_1$.
\end{cor}

\begin{proof}[Proof of Theorem\;\ref{thm:4-8}]
The norm $\NORM\cdot\NORM$ is constructed as in the proof of
Theorem\;\ref{thm:4-5} where we begin with $X =\langle b_0\rangle
\oplus_\infty [(b_i)]_{i\in\omega\setminus \{0\}}$ and
$(b_i)_{i\in\omega}$ is bimonotone. Everything we did in the proof of
Theorem\;\ref{thm:4-5} remains valid and in addition we have the use
of \eqref{eqn:4-11b}. It follows that not only do
$(x''_m)_{m\in\omega}$ and $(y''_n)_{n\in\omega}$ generate the same
type over $Y$, but so do $(x''_m)_{m\in\omega}$ and
$(-y''_n)_{n\in\omega}$ and thus, as in the case of
Theorem\;\ref{thm:4-5}, the proof reduces to the situation in
\cite{OS98b}. Hence, some subsequence of $(x''_m)_{m\in\omega}$ is an
$\ell_1$ basis.
\end{proof}

The arguments easily generalize to the case where $X$ is a subspace
of a space with a basis $(b_m)_{m\in\omega}$, and
$(e_n)_{n\in\omega}$ is generated by an array
$(x_i^n)_{n,i\in\omega}$, where for all $n,m$: $$\lim_{i\to\infty}
b_m^* (x_i^n) =0\,.$$

\section{Odds and Ends}\label{sec:odd}
\setcounter{equation}{0}
\setcounter{thm}{0}

In this section we first consider some stronger versions of
convergence one might hope for but, as we shall see, one cannot
always achieve. We also raise a number of open questions.

\subsection{Could We Get More?}

There are very many possible strengthenings of asymptotic models that
one could hope for. One such question is as follows:

Suppose we are
given a normalized basic sequence $(y_i)_{i\in\omega}$ and
$({\bfa}^i)_{i\in\omega}$.
{\sl Does there exist a subsequence
$(x_i)_{i\in\omega}$ of $(y_i)_{i\in\omega}$ with the following
property: for all $n\in\omega$, $(b_i)_{i\in n}\in [-1,1]^n$ and
$\eps >0$, there is an $N\in\omega$ so that if $N\le j_0<\ldots <j_{n-1}$,
$N\le k_0<\ldots
<k_{n-1}$ are integers and $Q\in\spart$, then}
$$\bigg| \Big\| \sum_{i\in n} b_i x\big( Q(j_i),{\bfa}^i\big) \Big\|
- \Big\| \sum_{i\in n} b_i x\big( Q(k_i),{\bfa}^i\big) \Big\| \bigg|
< \eps\ ?$$
Indeed, this is true if for each $i\in\omega$, ${\bfa}^i$
is finitely supported, for one can then take $(x_i)_{i\in\omega}$ to
be a subsequence of $(y_i)_{i\in\omega}$ generating a spreading model
$(\tilde{e}_i)_{i\in\omega}$.
The limit will exists in the above
sense (it will be just $\|\sum_{i\in n} b_i \tilde{f}_i \|$
where $(\tilde{f}_i)_{i\in\omega}$ is the normalized block
basis of $(\tilde{e}_i)_{i\in\omega}$ determined by the
${\bfa}^i$'s).

In general, however, this is false, even if $(y_i)_{i\in\omega}$ is
weakly null and ${\bfa}^i={\bfa}$ for all $i\in\omega$ and some
${\bfa}$. Indeed (cf. \cite[p.\,123]{LindenstraussTzafriri}) one can
embed $\ell_p\oplus\ell_2$ ($p\neq 2$) into a space $Y$ with a
normalized symmetric basis $(y_i)_{i\in\omega}$ in such a way that
the unit vector basis of $\ell_p\oplus\ell_2$ is equivalent to a
normalized block basis of the form $(y (P(i),{\bfa}))_{i\in\omega}$
where $|P(i)|\to\infty$ and ${\bfa}=(1,1,1,\ldots)$. Thus, for
appropriate $Q_1,Q_2\in\spart$ with $|Q_1(i)|, |Q_2(i)|\to\infty$,
every subsequence $(x_i)_{i\in\omega}$ of $(y_i)_{i\in\omega}$
contains block bases $\big(x (Q_1(i),{\bfa})\big)_{i\in\omega}$ and
$\big(x(Q_2(i), {\bfa})\big)_{i\in\omega}$ which are equivalent to
the unit vector basis of $\ell_p$ and $\ell_2$ respectively.

On the other hand, there are of course variations of our construction
of asymptotic models in Theorem\;\ref{thm:main} that do succeed. For
example, given a basic array $(x_i^n)_{n,i\in\omega}$, one might
stabilize $$\Big{\|}\sum_{i\in n} b_i x^{k(i,P)}
{\bfa}_{P(i)}^{k(i,P)} \Big{\|}$$ where the row now depends upon $i$
and $P\in\spart$. In this more general setting, one has that
$(e_i)_{i\in n}\in\{X_n\}$ iff there exists a block basic array
$(x_i^n)_{n,i\in\omega}$ and $k(i,P)$, ${\bfa}_{P(i)}^{j}$'s, so that
the above expression converges (as in Theorem\;\ref{thm:main}) to
$\|\sum_{i\in n} b_i e_i \|$. \jump \hspace*{\parindent} Indeed,
suppose for example that the tree $T_2=\{x_{(m_0,m_1)}: 0\le
m_0<m_1\}$ converges to $(e_1,e_2)$ as in  (\ref{rmk346}). Let
$x_i^0=x_{(i)}$, $x_i^1=x_{(0,i)}$ for $i>0$, $x_i^2=x_{(1,i)}$ for
$i>1$, and so on. (Notice that there is no need to define the first
part of each row.) Set $k(0,P):=0$ and $k(1,P):=j+1$ if $\min
P(0)=j$, and let ${\bfa}_{P(i)}^{j}:=(1,0,0,\ldots)$.

One could also relax the conditions defining a basic array $(x_i^n)$
by deleting the requirement that the rows be $K$-basic. This would
yield many more ``asymptotic models.'' For example every normalized
basic sequence $(x_i)$ in $X$ would be an ``asymptotic model'' of
$X$; take $(x_i^n) = (x_i)$ for all $n$. Proposition\;\ref{prp:37}
would also hold in this relaxed setting.

\subsection{Open Problems}\label{sec:problems}

\setcounter{thm}{\value{subsection}} \setcounter{subnummer}{0}
\beginsubnr{prb1}
$X$ is {\bf asymptotic $\ell_p$} (respectively, {\bf asymptotic
$c_0$}) if there exists $K$ so that for all $(e_i)_{i\in
n}\in\{X\}_n$, $(e_i)_{i\in n}$ is $K$-equivalent to the unit vector
basis of $\ell_p^n$ (respectively, $\ell_\infty^n$)
(see\;\cite{Maurey.etal}). Assume that there exists $K$ and $1\le
p\le\infty$ so that if $(e_i)_{i\in\omega}$ is an asymptotic model of
$X$, then $(e_i)_{i\in\omega}$ is $K$-equivalent to the unit vector
basis of $\ell_p$ ($c_0$, if $p=\infty$). {\sl Does $X$ contain an
asymptotic $\ell_p$ (or $c_0$) subspace?\/} The analogous problem for
spreading models is also open.
\end{subnr}
\beginsubnr{prb2}
Suppose $X$ has a basis and that there is a unique,
in the isometric sense, asymptotic model for all normalized block
basic arrays.
In this case, even if one replaces asymptotic model by
spreading model, it follows from Krivine's Theorem \cite{Krivine}
that this unique asymptotic model is $1$-equivalent to the unit
vector basis of $c_0$ or $\ell_p$ for some $1\le p<\infty$.
{\sl Must $X$ contain an isomorphic copy of this space?\/}
The analogous problem for spreading models is known to be true for the
case of $c_0$ and $\ell_1$
(see \cite{OS98b}).     
Also the asymptotic structure version of the question is true: if
$|\{X\}_2| =1$, then $X$ contains an isomorphic copy of
$c_0$ or $\ell_p$ (see\;\cite{Maurey.etal}).
\end{subnr}
\beginsubnr{prb3}
{\sl Can one stabilize the asymptotic models of a space $X$?\/}
Precisely, does there exist a basic sequence
$(x_i)_{i\in\omega}$ in $X$ so that for all block bases
$(y_i)_{i\in\omega}$ of $(x_i)_{i\in\omega}$, if
$(e_i)_{i\in\omega}$ is an asymptotic model of some normalized
block basic array of $(x_i)_{i\in\omega}$, then
$({e}_i)_{i\in\omega}$ is equivalent to an asymptotic model of a
normalized block basic array of $(y_i)_{i\in\omega}$\,?
We do not even know if there is some basic sequence $(x_i)_{i\in\omega}$ and
an asymptotic model $(e_i)_{i\in\omega}$ of $(x_i)_{i\in\omega}$
such that every block basis $(y_i)_{i\in\omega}$ of
$[(x_i)_{i\in\omega}]$ admits an asymptotic model
equivalent to $(e_i)_{i\in\omega}$.
The analogous questions for spreading models are open.
It is known that one can stabilize the
asymptotic structures $\{X\}_n$ for all $n\in\omega$ by passing to a
block basis (see\;\cite{Maurey.etal}).
\end{subnr}
\beginsubnr{prb4}
Assume that in $X$, every asymptotic
model $(e_i)_{i\in\omega}$ of any normalized basic block sequence
is $1$-unconditional (this is $\|\sum \pm a_i e_i\| = \|\sum a_i e_i\|$).
{\sl Does $X$ contain an unconditional basic sequence?
Does $X$ contain an asymptotically unconditional subspace?\/}
(i.e., a basic sequence
$(x_i)_{i\in\omega}$ so that for some $K<\infty$ and for all
$n\in\omega$, every block basis $(y_i)_{i\in n}$ of
$(x_i)_{i\in\omega\setminus n}$ is $K$-unconditional).
\end{subnr}
\beginsubnr{prb5}
{\sl For any space $X$, does there exist a finite chain of asymptotic
models $X=X_0,X_1,\ldots,X_n$, so that $X_{i+1}$ an asymptotic model
of $X_i$ (for $i\in n$) and $X_n$ is isomorphic to $c_0$ or $\ell_p$
for some $1\le p<\infty$\,?\/} The analogous problem for spreading
models is also open.
\end{subnr}

\beginsubnr{prb6}
For $1<p<\infty$, $\ell_p$ is arbitrarily distortable \cite{OS94}:
Given $K>1$ there exists an equivalent norm $\|\cdot\|$ on $\ell_p$
so that for all $X\subseteq \ell_p$, $(X,\|\cdot\|)$ is not
$K$-isomorphic to $\ell_p$. Is this true for asymptotic models as
well? {\sl Given $K>1$ (or for even some $K>1$) does there exist an
equivalent norm $\|\cdot\|$ on $\ell_p$ so that if $(e_i)_{i\in
\omega}$ is an asymptotic model of $(\ell_p,\|\cdot\|)$, then
$(e_i)_{i\in\omega}$ is not $K$-equivalent to the unit vector basis
of $\ell_p$\,?} The analogue for spreading models is also open.
\end{subnr}
\beginsubnr{prb7}
{\sl If $X$ has the property that every normalized bimonotone basic
sequence is an asymptotic model of $X$, does $X$ contain an
isomorphic copy of $c_0$\,?}
\end{subnr}

\end{document}